\documentclass[11pt,reqno]{amsart}

\usepackage[utf8]{inputenc}
\usepackage[margin=1in]{geometry}
\usepackage{amssymb}
\usepackage{mathtools}
\usepackage{pgfplots}
\pgfplotsset{compat=1.18}
\usetikzlibrary{positioning, arrows.meta}
\tikzset{
  ax/.style={-{Stealth[length=2mm]}},
  pt/.style={circle,fill,inner sep=0pt,minimum size=3.6pt},
  modcol/.style={green!45!black},
  modfill/.style={fill=green!12},
  slicecol/.style={blue!70!black},
  tent/.style={draw=blue!60!black,fill=blue!10},
  witness/.style={draw=blue!75!black,fill=blue!22,thick},
  rung/.style={gray!60},
}
\usepackage{eulervm,charter}
\usepackage{xcolor}
\definecolor{linkcol}{rgb}{0,0,0.55}
\definecolor{citecol}{rgb}{0,0.45,0}
\definecolor{urlcol}{rgb}{0.55,0,0}
\usepackage[colorlinks=true,linkcolor=linkcol,citecolor=citecol,urlcolor=urlcol]{hyperref}

\newcommand{\R}{\mathbb{R}}
\newcommand{\Z}{\mathbb{Z}}
\newcommand{\norm}[1]{\left\lVert #1 \right\rVert}
\newcommand{\db}{d_{B}}
\newcommand{\di}{d_{\mathrm{I}}}
\newcommand{\D}[1]{\mathcal{D}_{#1}}
\newcommand{\Dh}[1]{\widehat{\mathcal{D}}_{#1}}
\DeclareMathOperator{\diam}{diam}
\newcommand{\eps}{\varepsilon}
\newcommand{\LC}{\Lambda}

\newcommand{\Mnm}{\mathcal{M}_2^{n,m}}
\newcommand{\dM}{d_{\mathcal{M}}}
\newcommand{\Rk}{\operatorname{Rk}}
\newcommand{\asdim}{\operatorname{asdim}}

\newcommand{\RN}{\mathcal{R}_N}
\newcommand{\one}{\mathbf{1}}

\newtheorem{theorem}{Theorem}[section]
\newtheorem{lemma}[theorem]{Lemma}
\newtheorem{proposition}[theorem]{Proposition}
\newtheorem{corollary}[theorem]{Corollary}
\theoremstyle{definition}
\newtheorem{definition}[theorem]{Definition}

\theoremstyle{remark}
\newtheorem{remark}[theorem]{Remark}

\title[Coarse Embeddings Beyond Doubling]{Coarse Embeddings Beyond Doubling:\\ Persistence Diagrams and Multiparameter Modules}

\author{Sushovan Majhi}
\address{Data Science Program, George Washington University, Washington D.C., USA}
\email{s.majhi@gwu.edu}

\author{Atish Mitra}
\address{Department of Mathematical Sciences, Montana Technological University, Butte, MT, USA}
\email{amitra@mtech.edu}

\author{\v{Z}iga Virk}
\address{Faculty of Computer and Information Science, University of Ljubljana, Slovenia}
\email{zigavirk@gmail.com}


\begin{document}

\begin{abstract}
We study coarse embeddings into Hilbert space and the asymptotic dimension of spaces of persistence diagrams and multiparameter persistence modules.
We build coarse embeddings from admissible tent covers of bounded multiplicity at every scale, which doubling spaces carry but which do not require doubling: stacking them across dyadic scales gives an explicit $1$-Lipschitz coarse embedding, with lower gauge $c\sqrt{\tau t}$, or nearly linear for slowly decaying weights, and the same covers bound the asymptotic dimension by the multiplicity.
The spaces of persistence diagrams are not doubling but carry such covers, and for diagrams with at most $n$ points our covers improve the constant of the known landmark embedding by a factor exponential in $n$.
Slice by slice, the construction embeds multiparameter modules with boundedly many bars per slice under the matching distance on any finite net of slices, a class that contains the modules with at most $n$ generators.
Under the continuous distances, bounding bars or persistent Betti numbers does not suffice: two-parameter modules with one bar on every slice have infinite asymptotic dimension, and neither they nor the modules of bounded support and bounded rank admit a coarse embedding into Hilbert space, under the matching or the interleaving distance.
\end{abstract}

\maketitle

\section{Introduction}\label{sec:intro}

Persistent homology turns a filtration into a persistence module, and in one parameter into a persistence diagram; these invariants form metric spaces.
Many learning pipelines that use these invariants first map them into a Hilbert space, and whether such a map can keep far-apart inputs far apart is a coarse-geometric question: does the space coarsely embed, and what is its asymptotic dimension?
For diagrams with at most $n$ points, Mitra and Virk proved that it does, with asymptotic dimension $2n$~\cite{Mitra2021,MitraVirk2024corrigendum}.

One route to an explicit embedding is Assouad's theorem, which embeds snowflakes of doubling spaces bi-Lipschitzly into Euclidean space~\cite{Assouad1983,Heinonen2001,NaorNeiman2012}.
Diagram spaces are not doubling: the bottleneck distance collapses the diagonal to a point, so a ball about it contains infinitely many diagrams at mutual distance at least half its radius.
What they do have is admissible covers of bounded multiplicity at every scale, and the first aim of this paper is an explicit embedding, with a lower gauge, built from such covers alone.

With two parameters, a module no longer has a barcode, and two distances compete: the interleaving distance $\di$, over prime fields the largest stable one~\cite{Lesnick2015} but NP-hard to compute~\cite{BjerkevikBotnanKerber2020}, and the matching distance $\dM$, which compares the barcodes of one-parameter slices and can be computed~\cite{CerriDiFabioFerriEtAl2013,KerberLesnickOudot2019}.
The question is then whether a bound on a module's complexity implies embeddability, and the second aim is to address it for bars per slice, presentation size, and persistent Betti numbers.
The interleaving distance on modules of bounded presentation is treated in a companion paper~\cite{MajhiMitraVirk-asdim}.

\subsection*{Contribution}

The paper has a metric-space half and a persistence half.
The metric-space half dispenses with doubling: our construction needs only covers by \emph{tents}---distance-to-landmark bumps---that are admissible and of bounded multiplicity at every dyadic scale.
Stacking such covers across scales gives, on any metric space that carries them, a $1$-Lipschitz coarse embedding into a Hilbert space, $\ell^2$ when the space is separable, with an explicit, growing lower gauge $\rho_-(t)=c\sqrt{\tau t}$ (Theorem~\ref{thm:multi}\,(i)), nearly linear once the scale weights decay more slowly (Theorem~\ref{thm:multi}\,(ii)).
Doubling spaces carry such covers (Proposition~\ref{prop:net}), so the doubling case is included, but doubling is not needed.

The construction extends that of Mitra and Virk, which uses the same tents on a fixed lattice with another schedule of scales (Remark~\ref{rem:mv}), and explains why theirs works: $\D n$ carries admissible covers of bounded multiplicity at every scale, because points of persistence below the scale are absorbed by the diagonal (Lemma~\ref{lem:config}).
Since only admissibility is required, the landmarks may be placed where the data concentrate.

The persistence half asks what this buys for modules, and the answer depends on how their complexity is bounded.
\begin{enumerate}
\item \emph{Diagrams and sliced modules embed quantitatively.}
The diagrams $\D n$ embed with an explicit gauge, from covers (Lemma~\ref{lem:config}) that improve the gauge constant of the known landmark embedding by a factor exponential in $n$ (Theorem~\ref{thm:diagrams}, Remark~\ref{rem:mv}); by stacking slices, so do the modules with boundedly many bars per slice under the matching distance restricted to any finite net of slices (Theorem~\ref{thm:modules}).
\item \emph{Generators bound bars, but bars alone do not suffice.}
A module with at most $n$ generators has at most $n$ bars on every slice, whatever its relations (Lemma~\ref{lem:slice-cardinality}), so the sliced embeddings apply to bounded presentations.
Under the continuous distances, bounding bars instead fails: the modules with a \emph{single} bar on every slice contain isometric $\ell^\infty$-cubes of every dimension (Theorem~\ref{thm:staircase}), so under the matching and the interleaving distance alike their asymptotic dimension is infinite; since these cubes together contain every finite metric space isometrically, these modules admit no coarse embedding into Hilbert space either (Corollary~\ref{cor:staircase-nonembed}), nor into any Banach space of nontrivial type (Remark~\ref{rem:banach}).
Bounding bars per slice, the natural guess from one parameter, does not control the coarse geometry in two.
\item \emph{Rank bounds are too weak.}
Bounding the persistent Betti numbers, the invariant that does generalize from one parameter, gives no coarse embedding into Hilbert space: not in one parameter, where the bound allows any number of bars as long as few are alive at a time, and not for two-parameter modules of bounded support, under either distance (Theorems~\ref{thm:nonembed-diagrams} and~\ref{thm:nonembed-modules}); the obstruction is that of Dranishnikov, Gong, Lafforgue and Yu~\cite{DGLY2002}.
\end{enumerate}
On finite nets, bar bounds and presentation bounds alike yield quantitative embeddings and finite asymptotic dimension (Theorem~\ref{thm:modules}); under the continuous distances, the spaces cut out by bar bounds and by rank bounds admit no coarse embedding.
Whether presentation bounds give finite asymptotic dimension under the continuous matching distance is left open; the companion paper~\cite{MajhiMitraVirk-asdim} treats it together with the interleaving distance.

\subsection*{Related work}
\emph{Embeddings beyond doubling.}
Beyond doubling, a natural class is that of finite Assouad--Nagata dimension, which contains the doubling spaces~\cite{LangSchlichenmaier2005}, has a Hurewicz-type theory~\cite{BDLM2008}, and, in its asymptotic form, gives a Hilbert compression exponent of $1$~\cite{Gal2008}; finite asymptotic dimension alone gives a coarse embedding into Hilbert space~\cite[Ex.~11.5]{Roe2003},~\cite{BellDranishnikov2008}.
A multi-scale admissible family of bounded multiplicity witnesses finite asymptotic Assouad--Nagata dimension (Proposition~\ref{prop:covers-asdim}), so that the existence of coarse embeddings for the spaces treated here is consistent with that theory.
The tent construction adds the map itself, evaluated with a distance oracle on a given cover, with explicit constants, and, through the covers, the reason why the non-doubling diagram spaces embed.

\emph{Spaces of persistence diagrams.}
Mitra and Virk proved that $\D n$ coarsely embeds into Hilbert space with asymptotic dimension $2n$~\cite{Mitra2021,MitraVirk2024corrigendum}, and later gave explicit landmark embeddings~\cite{Mitra2024}. The space of finite diagrams under the bottleneck distance admits no coarse embedding~\cite{BubenikWagner2020,Mitra2021} (the space of all diagrams contains every bounded separable metric space isometrically~\cite{BubenikWagner2020}), nor do the diagram spaces with the $p$-Wasserstein distance for $2<p<\infty$~\cite{Wagner2021}; Bell, Lawson, Pritchard and Yasaki show that the space of finite diagrams fails Yu's property~A under the $p$-Wasserstein distances, $p<\infty$~\cite{BLPY2019}.
At bounded cardinality the finer, bi-Lipschitz question has been studied as well: Carri\`ere and Bauer showed that diagrams admit no bi-Lipschitz embedding into a finite-dimensional Euclidean space and that the distortion of several standard embeddings into Hilbert space grows with the cardinality~\cite{Carriere2019-nx}, and Bate and Garcia-Pulido embed unordered tuples bi-Lipschitzly into Hilbert space under the partial transportation metric~\cite{BateGarciaPulido2024}, which on diagrams with at most $n$ points is bi-Lipschitz to the bottleneck distance with constants depending on $n$.
Zava charts the same coarse-versus-bi-Lipschitz divide for subspaces of the Gromov--Hausdorff space~\cite{Zava2025}.
So bi-Lipschitz embeddings of bounded-cardinality diagrams into Hilbert space are known.
The landmark embeddings of Mitra and Virk~\cite{Mitra2024} and the tent embeddings below trade that control for maps evaluated with a distance oracle once the landmarks are chosen, each coordinate a function of the distance to a landmark.
The tent maps are coarse embeddings, nearly linear with the weights of Theorem~\ref{thm:multi}\,(ii), and not bi-Lipschitz on unbounded spaces (Remark~\ref{rem:distortion}); in exchange the tent construction applies to any space with admissible covers of bounded multiplicity at every scale, doubling or not, and extends slice by slice to multiparameter modules, beyond the reach of constructions specific to tuples.

\emph{Relation to Mitra--Virk.}
Most results of Mitra and Virk on diagram spaces have counterparts here, generalized or sharpened; Table~\ref{tab:mv} lists them.
Two are not recovered: the exact value $\asdim\D n=2n$ for $n\ge2$, because the multiplicities of $n$-tuple covers multiply, and the injective embedding of diagrams on a bounded domain into $\R^{n(n+1)}$~\cite[Ex.~5.4]{Mitra2024}, a different construction.

\begin{table}[t]
\centering\small
\renewcommand{\arraystretch}{1.15}
\begin{tabular}{p{0.40\textwidth}p{0.53\textwidth}}
\hline
Mitra--Virk & This paper \\
\hline
$\D n$ coarsely embeds into Hilbert space, for $\db$ and $W_p$ \cite[Thms.~2.12 and~3.2]{Mitra2021} &
explicit $1$-Lipschitz coarse embeddings of every space with admissible covers of bounded multiplicity (Thm.~\ref{thm:multi}), $\D n$ among them, for $\db$ and $W_p$ (Thm.~\ref{thm:diagrams}, Rem.~\ref{rem:mv}) \\
$\asdim\D n=2n$ \cite[Thm.~3.2]{Mitra2021}, \cite{MitraVirk2024corrigendum} &
$\asdim\D n\le3^n-1$ from the same covers, exact for $n=1$ (Prop.~\ref{prop:covers-asdim}, Lem.~\ref{lem:config}) \\
explicit coarse and uniform embeddings of $\D n$ with step gauges \cite[Thm.~4.3]{Mitra2024} &
the same construction for any schedule and any admissible cover (Rem.~\ref{rem:schedules}), cover constant larger by $3(2/\sqrt3)^n$, nearly linear gauge (Thm.~\ref{thm:multi}\,(ii), Rem.~\ref{rem:mv}) \\
bounded domains, Euclidean targets \cite[Thm.~5.1]{Mitra2024} &
bounded frames (Rem.~\ref{rem:complex}) \\
finite diagrams $\D{<\infty}$: infinite asymptotic dimension \cite[Cor.~2.11]{Mitra2021} &
two-parameter modules with one bar on every slice, under $\di$ and $\dM$ (Thm.~\ref{thm:staircase}) \\
$\D{<\infty}$: no coarse embedding into Hilbert space \cite[Thm.~4.3]{Mitra2021}, also \cite{BubenikWagner2020} &
diagrams with at most $n$ bars alive at a time (Thm.~\ref{thm:nonembed-diagrams}); two-parameter modules of bounded support with rank at most $n$, under $\di$ and $\dM$ (Thm.~\ref{thm:nonembed-modules}); modules with one bar on every slice, also for Banach targets of nontrivial type (Cor.~\ref{cor:staircase-nonembed}, Rem.~\ref{rem:banach}) \\
--- & sliced multiparameter modules (Thm.~\ref{thm:modules}) \\
\hline
\end{tabular}
\caption{Results of Mitra and Virk on diagram spaces and their counterparts here.}
\label{tab:mv}
\end{table}

\emph{Multiparameter persistence and learning.}
A large literature vectorizes multiparameter modules for learning: multiparameter landscapes~\cite{Vipond2020,Vipond2021}, images~\cite{Carriere2020-mpi}, signed measures and barcodes~\cite{Botnan2022signed,Loiseaux2023}, the differentiable family based on generalized rank invariant landscapes, \textsc{Gril}~\cite{Xin2023GRIL,Mukherjee2024}, and learned vectorizations of cubical multiparameter persistence~\cite{Korkmaz2025-cumperlay}.
These features are stable in one direction, their distances bounded above by the distance between modules; the two-sided question asked here complements them.
On the applied side, the closed-form landmark embeddings PLACE~\cite{PaperI}, on uniform grids, and their data-adaptive successor PALACE~\cite{PaperII} build one-parameter classifiers with a certified lower gauge, and COMPLEX~\cite{Majhi2026-complex} lifts them to multiparameter modules by stacking slices over a fixed near-diagonal net, with a lower gauge under a checkable coherence hypothesis. The tent covers below are the coarse-geometric counterpart of their landmark placement (Remarks~\ref{rem:mv} and~\ref{rem:complex}).

\subsection*{Organization}

Section~\ref{sec:prelim} fixes notation, including the matching distance; Section~\ref{sec:tents} develops the tent-cover embedding, for doubling spaces and beyond; Section~\ref{sec:persistence} applies it to diagrams and to modules under the matching distance on a finite net of slices.
Section~\ref{sec:slices} introduces modules of bounded presentation, proves the slice-cardinality lemma, and shows, via staircases, that one bar per slice is compatible with infinite asymptotic dimension and with the failure of coarse embeddability.
Section~\ref{sec:rank} proves the rank-condition obstruction to coarse embedding, and Section~\ref{sec:disc} discusses the open question.

\section{Preliminaries}\label{sec:prelim}

We fix the persistence and coarse-geometry conventions used below; all are standard.
The letter $d$ denotes the number of parameters, and also a metric or the death time of a bar where no confusion can arise.

\subsection{Multiparameter modules}\label{sec:modules}

Fix an integer $d\ge1$ and a field $k$, and order $\R^d$ componentwise: $x\le y$ when $x_i\le y_i$ for every $i$.
A \emph{$d$-parameter persistence module} is a functor $M\colon(\R^d,\le)\to\mathrm{Vec}_k$: a vector space $M(x)$ at every \emph{grade} $x\in\R^d$ and a linear \emph{structure map} $M(x\le y)\colon M(x)\to M(y)$ for every comparable pair, with $M(x\le x)=\mathrm{id}$ and $M(y\le z)\circ M(x\le y)=M(x\le z)$.
A \emph{morphism} $f\colon M\to N$ is a natural transformation, a family of linear maps $f_x\colon M(x)\to N(x)$ commuting with the structure maps; modules and morphisms form an abelian category, and kernels, images and cokernels are computed gradewise.
Sections~\ref{sec:prelim}--\ref{sec:persistence} are written for general $d$, except where a statement is marked as two-parameter; from Section~\ref{sec:slices} on, $d=2$.

\emph{Three examples recur throughout.}
The \emph{free module of rank one} at $\gamma\in\R^d$ is $k(-\gamma)$, equal to $k$ at every $x\ge\gamma$ and $0$ elsewhere, with identity structure maps where both ends are nonzero; its support is the closed orthant $\gamma+\R^d_{\ge0}$.
A \emph{bar} is a quotient $k(-\gamma)/k(-\gamma')$ with $\gamma\le\gamma'$, equal to $k$ on $\{x\ge\gamma,\ x\not\ge\gamma'\}$ and $0$ elsewhere; in one parameter this is the interval module on $[\gamma,\gamma')$.
More generally, for an \emph{upset} $U\subseteq\R^d$ (a set with $x\in U$, $y\ge x\Rightarrow y\in U$) the \emph{indicator module} $k[U]$ is $k$ on $U$ and $0$ elsewhere, with identity structure maps on $U$; in two parameters, when $U$ is generated by finitely many points, it is the staircase module of Section~\ref{sec:staircase}.
For $\eps\ge0$ the \emph{shift} $M[\eps]$ is the module $x\mapsto M(x+\eps\one)$, $\one=(1,\dots,1)$; the structure maps $M(x\le x+\eps\one)$ assemble into a morphism $M\to M[\eps]$, and $k[U][\eps]=k[U-\eps\one]$.

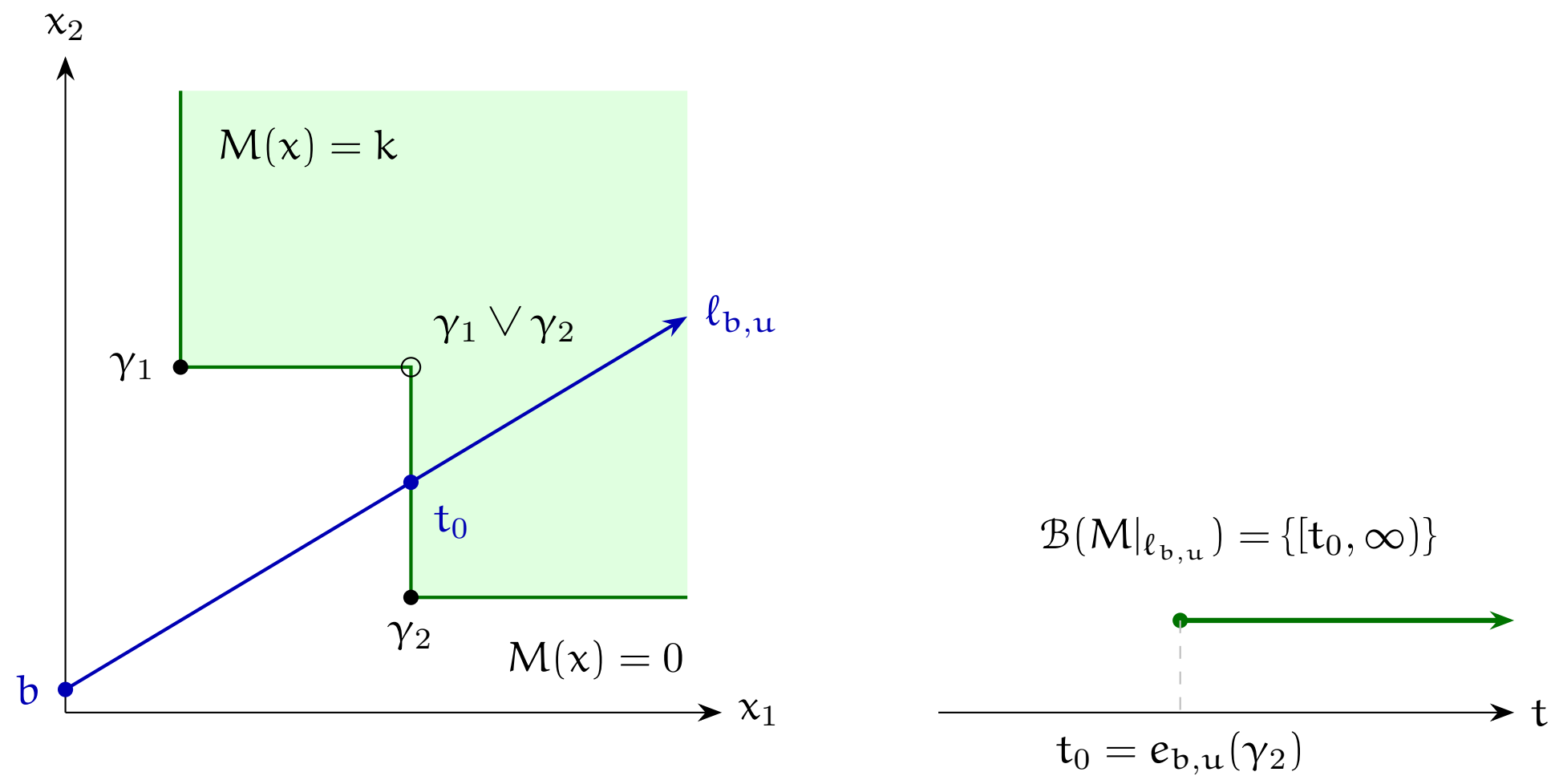
\begin{figure}[t]
\centering
\begin{tikzpicture}[x=0.95cm,y=0.95cm,font=\small]
  \begin{scope}
    \path[modfill] (1,3) -- (1,5.4) -- (5.4,5.4) -- (5.4,1) -- (3,1) -- (3,3) -- cycle;
    \draw[ax] (0,0) -- (5.7,0) node[right] {$x_1$};
    \draw[ax] (0,0) -- (0,5.7) node[above] {$x_2$};
    \draw[modcol,thick] (1,5.4) -- (1,3) -- (3,3) -- (3,1) -- (5.4,1);
    \node[pt,label=left:{$\gamma_1$}] at (1,3) {};
    \node[pt,label=below:{$\gamma_2$}] at (3,1) {};
    \node[circle,draw,inner sep=0pt,minimum size=4.4pt,label=above right:{$\gamma_1\vee\gamma_2$}] at (3,3) {};
    \node at (2.1,4.9) {$M(x)=k$};
    \node at (4.6,0.45) {$M(x)=0$};
    \draw[slicecol,thick,ax] (0,0.2) -- (5.4,3.44) node[right] {$\ell_{b,u}$};
    \node[pt,slicecol,label={[slicecol]left:$b$}] at (0,0.2) {};
    \node[pt,slicecol,label={[slicecol]below right:$t_0$}] at (3,2.0) {};
  \end{scope}
  \begin{scope}[xshift=7.2cm]
    \draw[ax] (0,0) -- (5,0) node[right] {$t$};
    \draw[modcol,very thick,ax] (2.1,0.8) -- (5,0.8);
    \node[pt,modcol] at (2.1,0.8) {};
    \draw[gray!50,dashed] (2.1,0.8) -- (2.1,0);
    \node at (2.1,-0.35) {$t_0=e_{b,u}(\gamma_2)$};
    \node at (2.6,1.5) {$\mathcal B(M|_{\ell_{b,u}})=\{[t_0,\infty)\}$};
  \end{scope}
\end{tikzpicture}
\caption{A two-parameter module and one of its slices. Left: the indicator module $M=k[U]$ of the upset generated by $\gamma_1=(1,3)$ and $\gamma_2=(3,1)$, with two generators and one relation at the join $\gamma_1\vee\gamma_2=(3,3)$, where the two generators become dependent; the shaded region is the support. The slice $\ell_{b,u}$ enters the support at the entry time $t_0=e_{b,u}(\gamma_2)$ of \eqref{eq:entry}, where it crosses the staircase. Right: the slice barcode, a single essential bar, as on every slice. Staircase modules with $q$ corners (Section~\ref{sec:staircase}) look the same with $q$ steps.}
\label{fig:module}
\end{figure}

\emph{Presentations.}
Following \cite[Def.~7.2]{BotnanLesnick2023}, a \emph{free presentation} of $M$ is a morphism $\partial\colon F_1\to F_0$ between direct sums of free modules of rank one with $\operatorname{coker}\partial\cong M$; the grades of the summands of $F_0$ are the \emph{generators} of the presentation and those of $F_1$ its \emph{relations}.
A module $M$ is \emph{finitely presented} if such a $\partial$ exists with $F_0,F_1$ finitely generated, and we assume this throughout; it implies that $M$ is pointwise finite-dimensional (p.f.d.).
By \cite[Thm.~7.6 and Rmk.~7.7(iii)]{BotnanLesnick2023} every finitely presented $M$ admits a minimal free resolution $0\to F_d\to\dots\to F_1\to F_0\to M\to0$, and in particular a minimal presentation $\partial_1\colon F_1\to F_0$, both unique up to isomorphism; we write $\beta_i^M$ for the multiset of grades of a basis of $F_i$ in the resolution, so that $|\beta_i^M|=\operatorname{rank}F_i$.
The generator and relation grades $\beta_0^M\cup\beta_1^M$, where the Betti numbers $\beta_0$ and $\beta_1$ are supported, are the \emph{Betti positions} of $M$.
Concretely, $x$ occurs in $\beta_0^M$ with multiplicity $\dim\bigl(M(x)/\sum_{y\le x,\,y\ne x}\operatorname{im}(M(y)\to M(x))\bigr)$, a free cover $F_0\to M$ is minimal exactly when its basis grades are $\beta_0^M$, and then $\beta_1^M=\beta_0^Z$ for the syzygy module $Z=\ker(F_0\to M)$, a minimal presentation being a minimal cover of $M$ followed by one of $Z$; these facts are standard (see~\cite[\S\S4.2 and~7]{BotnanLesnick2023} and~\cite{MillerSturmfels2005}).
In Figure~\ref{fig:module} the module $k[U]$ for the upset generated by $\gamma_1,\gamma_2$ has $\beta_0=\{\gamma_1,\gamma_2\}$ and $\beta_1=\{\gamma_1\vee\gamma_2\}$, the join being the componentwise maximum: the two generators become dependent exactly where both are present.

\emph{Slices and the fibered barcode.}
A \emph{slice} is a line $\ell\subset\R^d$ of positive slope; the product order restricts to a total order on $\ell$, so $M|_\ell\colon\ell\to\mathrm{Vec}_k$ is a one-parameter module.
It is convenient to parametrize a slice once and for all: a line of positive slope meets the orthant $\R^d_{\ge0}$ first at a point $b$ with $\min_ib_i=0$, and $\ell_{b,u}=\{b+tu:t\in\R\}$ where $u$ is a unit vector with positive coordinates, names it exactly once, $t$ being the signed Euclidean arclength from $b$.
A free module $k(-\gamma)$ restricts to the ray $[e_{b,u}(\gamma),\infty)$ with \emph{entry time}
\begin{equation}\label{eq:entry}
  e_{b,u}(\gamma)=\min\{t\in\R:b+tu\ge\gamma\}=\max_{1\le i\le d}\,(\gamma_i-b_i)/u_i ,
\end{equation}
which is the one computation behind every slice statement below: births on a slice are entry times of generators, and deaths are entry times of relations, since restricting a presentation to a slice gives a presentation of the restriction.
By the structure theorem for one-parameter modules \cite[Thm.~3.4]{BotnanLesnick2023}, $M|_\ell$ is a finite direct sum of interval modules, and the multiset of intervals is its \emph{barcode} $\mathcal B(M|_\ell)$; the assignment $\ell\mapsto\mathcal B(M|_\ell)$ is the \emph{fibered barcode} of $M$.
It determines the \emph{rank invariant} $\Rk M(x,y)=\operatorname{rank}M(x\le y)$, $x\le y$, since a comparable pair with $x<y$ in every coordinate lies on a slice and the rank there is the number of bars covering both parameters; the converse, and the degenerate pairs that lie on no slice, are treated in two parameters in the companion paper~\cite{MajhiMitraVirk-asdim}.
In one parameter, the barcode is a complete invariant of the isomorphism type; in two or more, neither the fibered barcode nor the rank invariant is~\cite{CarlssonZomorodian2009-multidim}.
We write $\D n$ for the space of persistence diagrams with at most $n$ off-diagonal points and $\db$ for the bottleneck distance; Definition~\ref{def:config} recalls both.

\subsection{Matching and interleaving distances}
An \emph{$\eps$-interleaving} of $M,N$ is a pair $f\colon M\to N[\eps]$, $g\colon N\to M[\eps]$ whose two composites are the internal $2\eps$-shift maps; the \emph{interleaving distance} is $\di(M,N)=\inf\{\eps\ge0: M,N\ \text{admit an}\ \eps\text{-interleaving}\}$, over prime fields the largest stable metric on modules~\cite{Lesnick2015} and NP-hard to compute~\cite{BjerkevikBotnanKerber2020}.
The matching distance instead aggregates the slices.
An $\eps$-interleaving shifts grades by $\eps\one$ and so moves each entry time~\eqref{eq:entry} on the slice $\ell_{b,u}$ by at most $\eps/w(u)$, where $w(u):=\min_iu_i$; this is the weight that makes each slice give a lower bound for the interleaving distance.

\begin{definition}[{\cite{CerriDiFabioFerriEtAl2013,Landi2018}}]\label{def:dmatch}
The \emph{matching distance} between p.f.d.\ modules $M,N$ is
\[
  \dM(M,N):=\sup_{b,u}\ w(u)\,\db\bigl(\mathcal B(M|_{\ell_{b,u}}),\mathcal B(N|_{\ell_{b,u}})\bigr),
\]
the supremum over all slices, with $\db$ the bottleneck distance.
Equivalently it is the unweighted supremum of $\db(\mathcal B(M\circ L),\mathcal B(N\circ L))$ over the \emph{normalized lines} $L(t)=b+tv$ with $b\in\R^d$, $v\in[1,\infty)^d$ and $\min_iv_i=1$; we write $\mathcal L$ for the set of these.
This is the form used by Botnan and Lesnick~\cite[Def.~6.11]{BotnanLesnick2023}, who let $v$ range over all of $[1,\infty)^d$: dividing $v$ by $\min_iv_i$ reparametrizes the same line and multiplies the bottleneck distance by $\min_iv_i$, so the supremum is unchanged. Moreover, $\dM\le\di$~\cite{Landi2018}.
\end{definition}

Either distance may be infinite between two modules, for instance when they differ in eventual dimension, so the spaces of modules below are extended metric spaces; the embeddings of Section~\ref{sec:tents} are then applied to each class of modules at finite mutual distance, with a base point in each.
For indicator modules of closed upsets, the interleaving distance is a geometric quantity: an $\eps$-interleaving between $k[U]$ and $k[V]$ exists exactly when $U+\eps\one\subseteq V$ and $V+\eps\one\subseteq U$ (Lemma~\ref{lem:upset-interleaving}, in two parameters), the computation behind Section~\ref{sec:slices}.
In one parameter both distances equal the bottleneck distance by the isometry theorem~\cite{Lesnick2015}; they differ in general as soon as $d\ge2$.

\subsection{Doubling and asymptotic dimension}
In this subsection and in Section~\ref{sec:tents}, $d$ denotes a metric.
Balls are open throughout, $B(x,r)=\{y:d(x,y)<r\}$, and $\overline B(x,r)$ is the closed ball.

\begin{definition}[Doubling metric space]\label{def:doubling}
$(X,d)$ is \emph{doubling} with \emph{doubling constant} $\mu\ge 1$ if every ball $B(x,r)$ is covered by at most $\mu$ balls of radius $r/2$.
Euclidean spaces, with any norm, are doubling.
\end{definition}

\begin{definition}[Asymptotic dimension]\label{def:asdim}
$(X,d)$ has \emph{asymptotic dimension} $\asdim X\le K$ if for every $R>0$ it admits a cover by uniformly bounded sets of \emph{$R$-multiplicity} at most $K+1$: every ball of radius $R$ meets at most $K+1$ members~\cite{BellDranishnikov2008}.
If the covers can be chosen with diameters at most $cR$ for $R$ at least some $R_0$, with $c$ independent of $R$, the \emph{asymptotic Assouad--Nagata dimension} is at most $K$~\cite{LangSchlichenmaier2005,BDLM2008}.
The definitions apply verbatim to extended metric spaces.
\end{definition}

Doubling is an all-scale condition, while asymptotic dimension is a coarse one; doubling implies finite asymptotic dimension, even finite Assouad--Nagata dimension~\cite{LangSchlichenmaier2005} (see also Proposition~\ref{prop:covers-asdim}), but not conversely, and the staircases in Section~\ref{sec:staircase} fail both.

\section{Tent-cover embeddings beyond doubling}\label{sec:tents}

We seek embeddings $\Phi$ of a metric space $(X,d)$ into Hilbert space that recover distance quantitatively from both sides: a pair of nondecreasing \emph{control functions}, or \emph{gauges}, $\rho_-\le\rho_+$ with
\[
  \rho_-\!\bigl(d(x,y)\bigr)\ \le\ \norm{\Phi(x)-\Phi(y)}\ \le\
  \rho_+\!\bigl(d(x,y)\bigr)\qquad(x,y\in X).
\]
The upper gauge is often Lipschitz; the content then lies in the lower gauge $\rho_-$, a certificate that separated points stay separated.
We take $\rho_-$ positive on a half-line $[\tau,\infty)$ and call $\Phi$ a \emph{coarse embedding} when, in addition, $\rho_-(t)\to\infty$; it is a \emph{uniform embedding} when $\rho_-$ is positive on all of $(0,\infty)$ and $\rho_+(t)\to0$ as $t\to0$, so that $\Phi$ and its inverse on the image are uniformly continuous.
The computational ingredient is a \emph{distance oracle}, a procedure that returns $d(x,y)$ for any two points, with no other access to them.

Fix a metric space $(X,d)$.
The embedding is assembled from elementary bumps, one per \emph{landmark}: to a center $p\in X$ and radius $r>0$ attach the \emph{tent}
\begin{equation}\label{eq:tent}
  \varphi_{p,r}(x):=\max\{\,r-d(p,x),\,0\,\},
\end{equation}
a pyramid of height $r$, supported on $\overline B(p,r)$ and $1$-Lipschitz in $d$.
A \emph{tent cover} $\LC=\{(p_i,r_i)\}_i$ is a cover of $X$ by such tents, which the embedding stacks, weighted by scale, in Theorem~\ref{thm:multi}. It can be chosen finite when $X$ is bounded and doubling, and it is countable whenever $X$ is separable and its multiplicity, defined below, is finite, since a countable dense set meets every ball $B(p_i,r_i)$ and each of its points lies in finitely many of them.
The lower bound turns on the cover's \emph{functional Lebesgue number}
\begin{equation}\label{eq:lebesgue}
  \lambda_0(\LC):=\inf_{x\in X}\ \sup_i\ \varphi_{p_i,r_i}(x),
\end{equation}
the least height the tents reach over $X$; when the multiplicity is finite the supremum is a maximum, and every point lies under some tent whose value there is at least $\lambda_0(\LC)$.
It plays the role of a Lebesgue number: since $\varphi_{p,r}(x)\ge h$ exactly when $x\in\overline B(p,r-h)$, for a tent cover of finite multiplicity one has $\lambda_0(\LC)\ge h$ precisely when the uniformly shrunk balls $\{\overline B(p_i,r_i-h)\}_i$ still cover $X$.
For a cover by tents of a common radius $r$ this reads $\lambda_0(\LC)=(r-\alpha)_+$, where $\alpha=\sup_x\inf_i d(x,p_i)$ is the covering radius of the landmark set, so a landmark set is described, for the purposes of this paper, by two numbers: its covering radius and its multiplicity.
A companion invariant, the cover's \emph{multiplicity}
\begin{equation}\label{eq:mult}
  \mathfrak m(\LC):=\sup_{x\in X}\ \#\{\,i:\varphi_{p_i,r_i}(x)>0\,\}
        =\sup_{x\in X}\ \#\{\,i:x\in B(p_i,r_i)\,\},
\end{equation}
counts the tents that are positive at a point; between two points at most $2\mathfrak m$ coordinates at a given scale change, so $\mathfrak m$ controls the Lipschitz constant of the stack.

To certify separation at a scale $\tau$, a cover must have tents no wider than $\tau$ that reach a height comparable to $\tau$ over every point; the next definition asks for both.

\begin{definition}[Admissibility]\label{def:admis}
Let $\tau>0$ and $\theta\in(0,1]$.
A tent cover $\LC$ is \emph{$(\tau,\theta)$-admissible} if every radius is $r_i\le\tau/2$ and its functional Lebesgue number is $\lambda_0(\LC)\ge\theta\tau/2$.
\end{definition}

\noindent The scale $\tau$ bounds the width of the tents, and the \emph{ratio} $\theta$ is the guaranteed height relative to that bound; rescaling $\tau$ rescales both conditions, so $\theta$ is the only shape parameter, and $\theta\le1$ since a tent of radius $r$ has height $r$.
Given two points at distance at least $\tau$, such a cover of finite multiplicity has over each a tent of value at least $\theta\tau/2$ there that vanishes at the other; these two coordinates are the entire lower-bound mechanism (see the proof of Theorem~\ref{thm:multi}).

Admissible tent covers exist in every metric space, produced by a separated \emph{net}; Theorem~\ref{thm:multi} needs bounded multiplicity, which doubling supplies (Proposition~\ref{prop:net}) and Lemma~\ref{lem:config} supplies on the non-doubling diagram spaces.

\begin{proposition}[Nets]\label{prop:net}
Let $(X,d)$ be a metric space, $\tau>0$ and $\theta\in(0,1)$.
A maximal $\eps$-separated set, $\eps=(1-\theta)\tau/2$ (pairwise distances $\ge\eps$), with uniform radii $r_i\equiv\tau/2$, is a $(\tau,\theta)$-admissible tent cover.
If $X$ is doubling with constant $\mu$, its multiplicity is at most $\mu^{k_\theta}$, $k_\theta=\lceil\log_2(2/(1-\theta))\rceil$; for $\theta=\tfrac12$ this is $\mu^2$.
\end{proposition}

\begin{proof}
A maximal $\eps$-separated set exists by Zorn's lemma, and by maximality every $x$ lies within distance less than $\eps$ of some center, under a tent of height greater than $\tau/2-\eps=\theta\tau/2$; thus $\lambda_0\ge\theta\tau/2$, which with radii $r_i\equiv\tau/2$ is $(\tau,\theta)$-admissibility.
For the multiplicity, the centers of the tents positive at a point $x$ lie in $B(x,\tau/2)$; $k_\theta$ applications of doubling cover that ball by $\mu^{k_\theta}$ balls of radius $\tau/2^{k_\theta+1}\le\eps/2$, any two points of each at distance less than $\eps$, so that each holds at most one center, and at most $\mu^{k_\theta}$ tents are positive at $x$.
\end{proof}

\noindent The ratio trades against the multiplicity: the gauge of Theorem~\ref{thm:multi} carries the factor $\theta/\sqrt{\mathfrak m}\ge\theta\mu^{-k_\theta/2}$, and $\theta=\tfrac12$ gives $1/(2\mu)$.

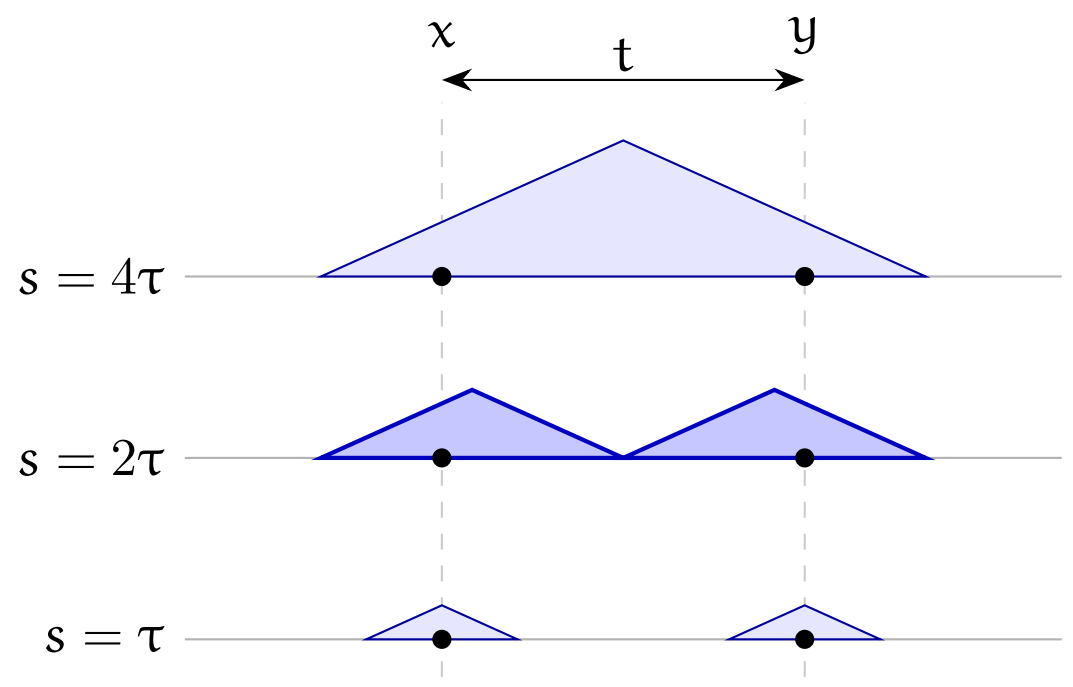
\begin{figure}[t]
\centering
\begin{tikzpicture}[x=1cm,y=1cm,font=\small]
  \def\px{2.3} \def\py{4.7}
  \node at (\px,4.0) {$x$};
  \node at (\py,4.0) {$y$};
  \draw[{Stealth[length=2mm]}-{Stealth[length=2mm]}] (\px,3.7) -- node[above,inner sep=1.5pt]{$t$} (\py,3.7);
  \draw[gray!40,dashed] (\px,-0.25) -- (\px,3.55);
  \draw[gray!40,dashed] (\py,-0.25) -- (\py,3.55);
  \begin{scope}[yshift=2.4cm]
    \draw[rung] (0.6,0) -- (6.4,0);
    \draw[tent] (1.5,0) -- (3.5,0.9) -- (5.5,0) -- cycle;
    \node[left] at (0.6,0) {$s=4\tau$};
    \node[pt] at (\px,0) {}; \node[pt] at (\py,0) {};
  \end{scope}
  \begin{scope}[yshift=1.2cm]
    \draw[rung] (0.6,0) -- (6.4,0);
    \draw[witness] (1.5,0) -- (2.5,0.45) -- (3.5,0) -- cycle;
    \draw[witness] (3.5,0) -- (4.5,0.45) -- (5.5,0) -- cycle;
    \node[left] at (0.6,0) {$s=2\tau$};
    \node[pt] at (\px,0) {}; \node[pt] at (\py,0) {};
  \end{scope}
  \begin{scope}[yshift=0cm]
    \draw[rung] (0.6,0) -- (6.4,0);
    \draw[tent] (1.8,0) -- (2.3,0.225) -- (2.8,0) -- cycle;
    \draw[tent] (4.2,0) -- (4.7,0.225) -- (5.2,0) -- cycle;
    \node[left] at (0.6,0) {$s=\tau$};
    \node[pt] at (\px,0) {}; \node[pt] at (\py,0) {};
  \end{scope}
\end{tikzpicture}
\caption{Three rungs of the dyadic ladder at a pair $(x,y)$ with $t=d(x,y)\in[2\tau,4\tau)$; tents on the rung $s$ have radius at most $s/2$, drawn to one vertical scale.
On the witnessing rung $s=2\tau$ (bold) a tent rises at least $\theta s/2$ over each point and vanishes at the other: these two coordinates carry the lower bound of Theorem~\ref{thm:multi}.
The coarse rung covers both points with one tent, and the fine rung separates them with shorter tents.}
\label{fig:ladder}
\end{figure}

A single admissible cover certifies separation only up to its own scale.
The \emph{dyadic ladder} removes that ceiling: on every rung $s_j=\tau2^j$, $j\ge0$, place an $(s_j,\theta)$-admissible cover of multiplicity at most $\mathfrak m$, and stack the covers with weights $2^{-j/2}$; we call the scales and the weights the \emph{schedule}.
The bands $[s_j,s_{j+1})$ partition the distances $[\tau,\infty)$, so a pair at distance $t\ge\tau$ lands in one band, $s_{j^\star}\le t<2s_{j^\star}$.
On that rung a tent rises at least $\theta s_{j^\star}/2$ over each point, and the other point, a distance $t\ge s_{j^\star}$ off, escapes it (Figure~\ref{fig:ladder}).
These two coordinates carry the lower bound, of order $2^{-j^\star/2}s_{j^\star}\approx\sqrt{\tau t}$.
Finer rungs separate the pair as well, with smaller contributions, and coarser rungs need not separate it at all.
The weights make the stack $1$-Lipschitz, and measuring every coordinate relative to a base point keeps it in $\ell^2$.
For doubling spaces such power-law gauges go back to Assouad~\cite{Assouad1983}; here the only hypothesis is admissible covers of bounded multiplicity, and once the covers are chosen the map is evaluated with a distance oracle.

\begin{definition}[Multi-scale admissible family]\label{def:msfamily}
Let $\tau>0$, $\theta\in(0,1]$ and $s_j=\tau2^j$, $j\ge0$. A \emph{multi-scale admissible family} of ratio $\theta$ and multiplicity $\le \mathfrak m$ on $(X,d)$ is a sequence $\mathcal F=\{\LC_j\}_{j\ge0}$ in which each $\LC_j$ is an $(s_j,\theta)$-admissible tent cover (Definition~\ref{def:admis}) of multiplicity $\le \mathfrak m$.
\end{definition}

\begin{theorem}[Multi-scale coarse embedding]\label{thm:multi}
Let $(X,d)$ be a metric space and $\mathcal F=\{\LC_j\}_{j\ge0}$ a multi-scale admissible family of ratio $\theta$ and multiplicity $\le \mathfrak m$ on $(X,d)$, at the scales $s_j=\tau2^j$.
Fix weights $a_j>0$ with $W:=\sum_{j\ge0}a_j^2<\infty$, a base point $x_0\in X$, and $\kappa:=(2\mathfrak m W)^{-1/2}$, and form the multi-scale embedding whose coordinates are tent values relative to $x_0$,
\begin{equation}\label{eq:ms}
  \Phi(x;\mathcal F):=\Bigl(\kappa\,a_j\bigl(\varphi_{p,r}(x)
  -\varphi_{p,r}(x_0)\bigr)\Bigr)_{\,j\ge0,\ (p,r)\in\LC_j}
\end{equation}
so that $\Phi(x_0)=0$.
Then $\Phi(\cdot;\mathcal F)$ is a $1$-Lipschitz map into $\ell^2(I)$, $I$ the disjoint union of the index sets of the $\LC_j$ (so into $\ell^2$ when $X$ is separable), and for $x,y\in X$ with $t=d(x,y)\ge\tau$,
\[
  \norm{\Phi(x;\mathcal F)-\Phi(y;\mathcal F)}\ \ge\ \tfrac{\theta}{\sqrt2}\,\kappa\,a_{j^\star}s_{j^\star},\qquad j^\star=\lfloor\log_2(t/\tau)\rfloor,\quad s_{j^\star}\in(t/2,t].
\]
In particular:
\begin{enumerate}
\item with the dyadic weights $a_j=2^{-j/2}$,
\begin{equation}\label{eq:snow}
  c_\theta\,\mathfrak m^{-1/2}\,\sqrt{\tau t}\ \le\
  \norm{\Phi(x;\mathcal F)-\Phi(y;\mathcal F)}_{\ell^2}\ \le\ \min\{t,\sqrt{\tau t}\},\qquad c_\theta=\tfrac{\theta}{4},
\end{equation}
where the bound $t$ holds for all $x,y$ and the other two for $t\ge\tau$;
\item with $a_j=(j+1)^{-(1+\delta)/2}$, $\delta>0$, the gauge is nearly linear: for $t\ge\tau$,
\[
  \norm{\Phi(x;\mathcal F)-\Phi(y;\mathcal F)}\ \ge\ \frac{c_{\theta,\delta}}{\sqrt{\mathfrak m}}\,\frac{t}{\bigl(1+\log_2(t/\tau)\bigr)^{(1+\delta)/2}},\qquad c_{\theta,\delta}=\frac{\theta}{4\sqrt{\zeta(1+\delta)}},
\]
where $\zeta(1+\delta)=\sum_{k\ge1}k^{-(1+\delta)}$ is a value of the Riemann zeta function, finite because $\delta>0$.
\end{enumerate}
\end{theorem}

\begin{proof}
Write $\eta_{j,p}:=\varphi_{p,r}(x)-\varphi_{p,r}(y)$; each tent is $1$-Lipschitz with values in $[0,s_j/2]$, so $|\eta_{j,p}|\le\min(t,s_j/2)\le t$, and at scale $j$ at most $2\mathfrak m$ of the $\eta_{j,p}$ are nonzero.
Hence
\[
  \norm{\Phi(x)-\Phi(y)}^2=\kappa^2\sum_{j\ge0}a_j^2\sum_p\eta_{j,p}^2
  \le\kappa^2\cdot 2\mathfrak m\,t^2\sum_{j\ge0}a_j^2=2\mathfrak m W\kappa^2t^2=t^2,
\]
and $\Phi$ is $1$-Lipschitz.
When $X$ is separable each $\LC_j$ is countable, as noted after \eqref{eq:tent}, so $I$ is countable.

For the lower bound, take $t\ge\tau$ and let $j^\star$ be the largest $j$ with $s_j\le t$, so $s_{j^\star}\in(t/2,t]$.
Write $s=s_{j^\star}$.
By $(s,\theta)$-admissibility some tent has $\varphi_{p,r}(x)\ge\theta s/2$, the maximum in \eqref{eq:lebesgue} being over the at most $\mathfrak m$ tents positive at $x$.
Then $d(p,x)\le r-\theta s/2\le(1-\theta)s/2$, so $d(p,y)\ge t-(1-\theta)s/2\ge(1+\theta)s/2>s/2\ge r$ and $\varphi_{p,r}(y)=0$.
In the same way some tent has $\varphi_{q,r'}(y)\ge\theta s/2$ and $\varphi_{q,r'}(x)=0$, and it is not the first, which vanishes at $y$.
These two coordinates give $\norm{\Phi(x)-\Phi(y)}^2\ge2\bigl(\kappa\,a_{j^\star}\theta s/2\bigr)^2$, which is the displayed bound.
For (i), $W=2$, $\kappa=1/(2\sqrt{\mathfrak m})$, and $\tfrac{\theta}{\sqrt2}\kappa\,2^{-j^\star/2}s_{j^\star}=\tfrac{\theta}{\sqrt2}\kappa\sqrt{\tau s_{j^\star}}\ge\frac{\theta\sqrt{\tau t}}{4\sqrt{\mathfrak m}}$, using $s_{j^\star}>t/2$.
For (ii), $W=\zeta(1+\delta)$, $s_{j^\star}>t/2$, and $j^\star\le\log_2(t/\tau)$ gives $a_{j^\star}\ge(1+\log_2(t/\tau))^{-(1+\delta)/2}$.

For the bound $\sqrt{\tau t}$ in (i), let $t\ge\tau$.
With $|\eta_{j,p}|\le\min(t,s_j/2)$ and $\kappa^2=1/(4\mathfrak m)$, the computation above gives $\norm{\Phi(x)-\Phi(y)}^2\le\tfrac12\sum_j2^{-j}\min(t,s_j/2)^2$.
The terms with $j\le j^\star+1$ have $s_j/2\le s_{j^\star}\le t$ and sum to $\sum_{j\le j^\star+1}\tau^22^j/8<\tau^22^{j^\star}/2=\tau s_{j^\star}/2\le\tau t/2$; the terms with $j\ge j^\star+2$ sum to at most $t^2\,2^{-(j^\star+2)}=t^2\tau/(4s_{j^\star})<\tau t/2$, since $t<2s_{j^\star}$.
Hence $\norm{\Phi(x)-\Phi(y)}^2<\tau t$.

The base point is needed: without it, $\Phi(x)$ has at every scale $j$ a coordinate at least $\theta\kappa a_js_j/2$, which is unbounded in $j$ for the weights of (i) and (ii).
With it, the computation above applied to $y=x_0$ gives $\norm{\Phi(x)}\le d(x,x_0)<\infty$.
\end{proof}

Unless stated otherwise, we use the dyadic weights of (i); with them every pair at distance at least $\tau$ is separated by at least $c_\theta\mathfrak m^{-1/2}\tau$, a per-scale floor.

\begin{remark}[Other scales; uniform embeddings]\label{rem:schedules}
The lower-bound argument uses only that the witnessing scale is at most $t$.
So Theorem~\ref{thm:multi} holds for any strictly increasing scales $(s_j)_{j\in J}$, $J=\mathbb N$ or $J=\Z$, with $(s_j,\theta)$-admissible covers of multiplicity $\le\mathfrak m$, weights with $W=\sum_{j\in J}a_j^2<\infty$ and $\kappa=(2\mathfrak mW)^{-1/2}$: $\Phi$ is $1$-Lipschitz, and $\norm{\Phi(x)-\Phi(y)}\ge\tfrac{\theta}{\sqrt2}\kappa\,a_{j^\star}s_{j^\star}$ with $s_{j^\star}$ the largest scale $\le t=d(x,y)$, whenever there is one.
With $J=\Z$, $s_j\to0$ as $j\to-\infty$ and $s_j\to\infty$ as $j\to\infty$ there is one for every $t>0$, and $\Phi$ is a uniform as well as a coarse embedding: for $s_j=\tau2^j$ and $a_j=2^{-|j|/2}$, $j\in\Z$, one has $W=3$ and
\[
  \norm{\Phi(x)-\Phi(y)}\ \ge\ \frac{\theta}{2\sqrt{6\mathfrak m}}\,\sqrt{\tau t}\quad(t\ge\tau),\qquad
  \norm{\Phi(x)-\Phi(y)}\ \ge\ \frac{\theta}{4\sqrt{6\mathfrak m\tau}}\,t^{3/2}\quad(t\le\tau),
\]
by the computation in the proof of Theorem~\ref{thm:multi}, with $a_{j^\star}s_{j^\star}=s_{j^\star}^{3/2}/\sqrt\tau$ for $j^\star<0$.
This is the setting of the coarse and uniform embeddings of \cite[Thm.~4.3]{Mitra2024}, built there from increasing scales and from scales decreasing to $0$ and combined into one map; the two-sided ladder needs covers at every small scale as well, which the diagram spaces have (Lemma~\ref{lem:config}).

On cones a single landmark set serves every scale: if $(X,d)$ has surjective dilations $\sigma_\lambda$ fixing $x_0$, with $d(\sigma_\lambda x,\sigma_\lambda y)=\lambda d(x,y)$, then $\varphi_{\sigma_\lambda p,\lambda r}\circ\sigma_\lambda=\lambda\varphi_{p,r}$, so rescaling one $(\tau,\theta)$-admissible cover by $\sigma_{2^j}$ gives a $(2^j\tau,\theta)$-admissible cover of the same multiplicity. Diagram and module spaces are cones, with the empty diagram and the zero module fixed; under a sliced distance the dilations must also preserve the slice net, as they preserve a fan of slices through the origin.
\end{remark}

\begin{remark}[Distortion]\label{rem:distortion}
With the dyadic weights and $t=d(x,y)\ge\tau$, Theorem~\ref{thm:multi}\,(i) bounds $\norm{\Phi(x)-\Phi(y)}$ between $c_\theta\mathfrak m^{-1/2}\sqrt{\tau t}$ and $\sqrt{\tau t}$: on pairs at distance at least $\tau$, $\Phi$ is a bi-Lipschitz embedding of the snowflake $\sqrt{\tau d}$, with distortion at most $4\sqrt{\mathfrak m}/\theta$.
With the weights of Theorem~\ref{thm:multi}\,(ii) the ratio of the gauges is only a power of $\log(t/\tau)$, so every space carrying a multi-scale admissible family of bounded multiplicity has Hilbert compression exponent $1$.
On an unbounded space neither map is bi-Lipschitz for $d$ itself, since $\norm{\Phi(x)-\Phi(y)}/d(x,y)\to0$ as $d(x,y)\to\infty$: the ratio is at most $\sqrt{\tau/t}$ in case (i), and the same computation bounds it by a constant multiple of $(1+\log(t/\tau))^{-\delta/2}$ in case (ii).
The $n$-point diagrams do admit bi-Lipschitz embeddings into Hilbert space~\cite{BateGarciaPulido2024}; the tent maps give up that optimum for evaluation with a distance oracle, and need only admissible covers.
\end{remark}

The covers that drive the embedding also bound the asymptotic dimension, so one family of covers gives both sides of the implication from finite asymptotic dimension to a coarse embedding~\cite[Ex.~11.5]{Roe2003},~\cite{BellDranishnikov2008}, with explicit constants.

\begin{proposition}[Admissible covers bound the asymptotic dimension]\label{prop:covers-asdim}
If $(X,d)$ carries a multi-scale admissible family of ratio $\theta$ and multiplicity $\le \mathfrak m$, then $\asdim X\le \mathfrak m-1$; in fact the covers can be chosen with diameter $<4R/\theta$ at every sufficiently large scale $R$, so the asymptotic Assouad--Nagata dimension is at most $\mathfrak m-1$.
In particular a doubling space of constant $\mu$ has $\asdim X\le\mu^2-1$, by the nets of Proposition~\ref{prop:net} with $\theta=\tfrac12$ at every scale.
\end{proposition}

\begin{proof}
Let $\tau$ be the base scale of the family, so that its scales are $s_j=\tau2^j$.
Given $R>0$, let $s=s_j$ be the least scale with $\theta s/2\ge R$, and let $\LC_j=\{(p_i,r_i)\}_i$ be the $(s,\theta)$-admissible cover of the family, with $r_i\le s/2$ and $\lambda_0(\LC_j)\ge\theta s/2\ge R$.
Put $C_i:=\{x:\varphi_{p_i,r_i}(x)\ge R\}=\overline B(p_i,r_i-R)$.
Since $\lambda_0\ge R$ and the supremum in \eqref{eq:lebesgue} is a maximum, the nonempty $C_i$ cover $X$; their diameters are at most $2r_i\le s$.
If the ball $B(x,R)$ meets $C_i$ at a point $y$, then $d(p_i,x)<(r_i-R)+R=r_i$, so the tent $\varphi_{p_i,r_i}$ is positive at $x$, which happens for at most $\mathfrak m$ indices $i$.
Thus $\{C_i\}_i$ is a uniformly bounded cover of $R$-multiplicity $\le\mathfrak m$, which is Definition~\ref{def:asdim} with $K+1=\mathfrak m$.
Since the scales double, the least $s_j\ge2R/\theta$ is less than $4R/\theta$ as soon as $2R/\theta\ge\tau$, which gives the diameter bound.
For a doubling space the family is given by the nets of Proposition~\ref{prop:net}, which exist at every scale.
\end{proof}

\section{Diagrams and sliced modules}\label{sec:persistence}

The cover theory of Section~\ref{sec:tents} never mentions persistence: the tents act on whatever metric space we supply, so what remains is to name the right space and to supply admissible covers of bounded multiplicity.
We do this twice: for $n$-point diagrams (Section~\ref{sec:diagrams}), then, by slicing, for multiparameter modules (Section~\ref{sec:multiparam}).

\subsection{\texorpdfstring{$n$}{n}-point diagrams}\label{sec:diagrams}

\begin{definition}[$n$-point diagrams]\label{def:config}
Following Mitra--Virk~\cite{Mitra2024}, the persistence diagrams on $n$ points are $\D{n}=\D{1}^{\,n}/S_n$, where $\D{1}=H\cup\{\Delta\}$ adjoins to the open half-plane $H=\{(b,d):b<d\}$ a single point $\Delta$ for the diagonal, with $\db(x,\Delta)=\tfrac12(d-b)$ for $x=(b,d)$ and $\db(x,y)=\min\{\norm{x-y}_\infty,\max(\db(x,\Delta),\db(y,\Delta))\}$.
The \emph{bottleneck distance} on $\D{n}$ is the induced $\ell^\infty$-matching $\db(A,B)=\min_{\sigma\in S_n}\max_i\db(a_i,b_{\sigma(i)})$, whose matchings may route points through $\Delta$.
\end{definition}

The diagonal makes $(\D{n},\db)$ non-doubling, already for $n=1$: $\db$ caps distances, $\db(x,y)\le\max\{\db(x,\Delta),\db(y,\Delta)\}$, so the ball $B(\Delta,r)$ runs unboundedly along the diagonal yet no finite family of $r/2$-balls covers it.
Nor does the failure need the diagonal itself: points of persistence between $r$ and $2r$ that are pairwise at least $r$ apart in $\ell^\infty$ are $r$-separated in $\db$ yet all within $2r$ of one another, for every $r>0$.
What saves the construction is that the cap is relative to the scale.
Theorem~\ref{thm:multi} asks only for admissible covers of bounded multiplicity, and $\D{n}$ supplies them. At scale $\tau$ the cap matters only for points of persistence below about $\tau$, which all lie within $\tau$ of the single point $\Delta$ and are absorbed by one tent centered at $\Delta$; about any point of persistence at least $\tau$, the $\db$-balls of radius at most $\tau$ are $\ell^\infty$-balls, with the doubling geometry of $\R^{2n}$.

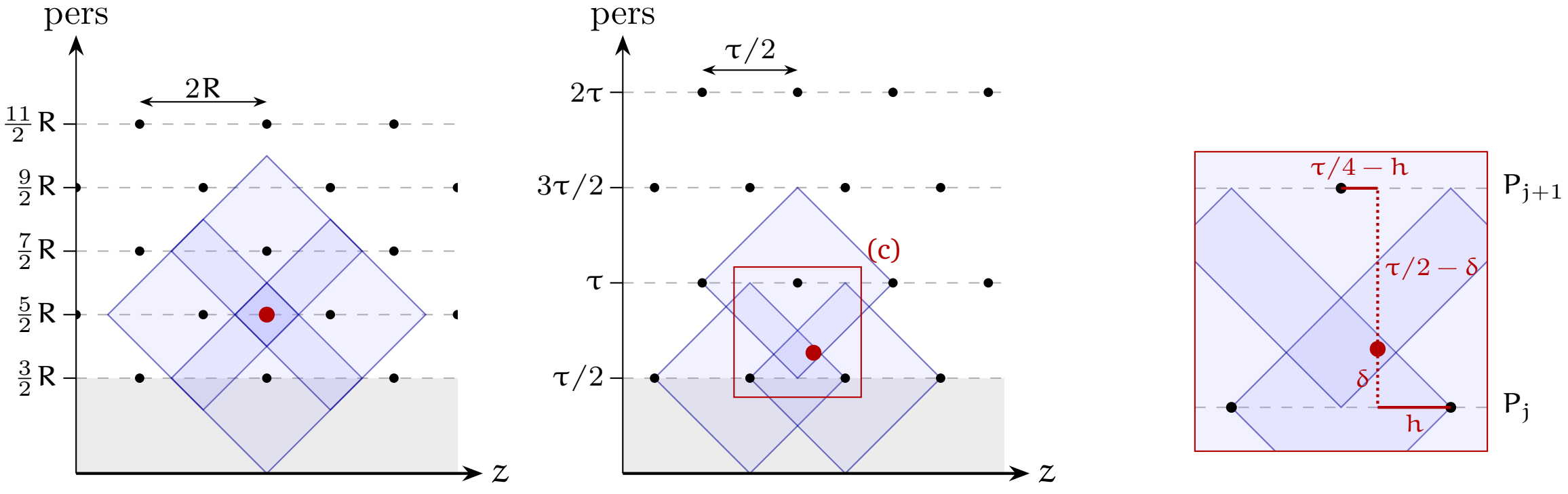
\begin{figure}[t]
\centering
\begin{tikzpicture}[x=0.62cm,y=0.62cm,font=\small]
  \begin{scope}[xshift=0.31cm]
    \begin{scope}
    \clip (-0.5,0) rectangle (5.5,6.3);
    \path[fill=gray!15] (-0.5,0) rectangle (5.5,1.5);
    \draw[rung,dashed] (-0.5,1.5) -- (5.5,1.5);
    \draw[rung,dashed] (-0.5,2.5) -- (5.5,2.5);
    \draw[rung,dashed] (-0.5,3.5) -- (5.5,3.5);
    \draw[rung,dashed] (-0.5,4.5) -- (5.5,4.5);
    \draw[rung,dashed] (-0.5,5.5) -- (5.5,5.5);
    \draw[tent,fill=blue!45,fill opacity=0.12,draw opacity=0.55] (0,2.5) -- (1.5,4) -- (3,2.5) -- (1.5,1) -- cycle;
    \draw[tent,fill=blue!45,fill opacity=0.12,draw opacity=0.55] (1,3.5) -- (2.5,5) -- (4,3.5) -- (2.5,2) -- cycle;
    \draw[tent,fill=blue!45,fill opacity=0.12,draw opacity=0.55] (1,1.5) -- (2.5,3) -- (4,1.5) -- (2.5,0) -- cycle;
    \draw[tent,fill=blue!45,fill opacity=0.12,draw opacity=0.55] (2,2.5) -- (3.5,4) -- (5,2.5) -- (3.5,1) -- cycle;
    \node[pt,minimum size=2.6pt] at (-0.5,4.5) {};
    \node[pt,minimum size=2.6pt] at (0.5,5.5) {};
    \node[pt,minimum size=2.6pt] at (-0.5,2.5) {};
    \node[pt,minimum size=2.6pt] at (0.5,3.5) {};
    \node[pt,minimum size=2.6pt] at (1.5,4.5) {};
    \node[pt,minimum size=2.6pt] at (2.5,5.5) {};
    \node[pt,minimum size=2.6pt] at (0.5,1.5) {};
    \node[pt,minimum size=2.6pt] at (1.5,2.5) {};
    \node[pt,minimum size=2.6pt] at (2.5,3.5) {};
    \node[pt,minimum size=2.6pt] at (3.5,4.5) {};
    \node[pt,minimum size=2.6pt] at (4.5,5.5) {};
    \node[pt,minimum size=2.6pt] at (2.5,1.5) {};
    \node[pt,minimum size=2.6pt] at (3.5,2.5) {};
    \node[pt,minimum size=2.6pt] at (4.5,3.5) {};
    \node[pt,minimum size=2.6pt] at (5.5,4.5) {};
    \node[pt,minimum size=2.6pt] at (4.5,1.5) {};
    \node[pt,minimum size=2.6pt] at (5.5,2.5) {};
    \end{scope}
    \draw[thick,-{Stealth[length=2mm]}] (-0.5,0) -- (5.9,0) node[right,inner sep=2pt] {$z$};
    \draw[ax] (-0.5,0) -- (-0.5,6.9) node[above,inner sep=1pt] {$\operatorname{pers}$};
    \draw (-0.5,1.5) -- (-0.7,1.5) node[left,inner sep=1.5pt,font=\scriptsize] {$\tfrac32R$};
    \draw (-0.5,2.5) -- (-0.7,2.5) node[left,inner sep=1.5pt,font=\scriptsize] {$\tfrac52R$};
    \draw (-0.5,3.5) -- (-0.7,3.5) node[left,inner sep=1.5pt,font=\scriptsize] {$\tfrac72R$};
    \draw (-0.5,4.5) -- (-0.7,4.5) node[left,inner sep=1.5pt,font=\scriptsize] {$\tfrac92R$};
    \draw (-0.5,5.5) -- (-0.7,5.5) node[left,inner sep=1.5pt,font=\scriptsize] {$\tfrac{11}2R$};
    \draw[{Stealth[length=1.3mm]}-{Stealth[length=1.3mm]}] (0.5,5.85) -- node[above,inner sep=1pt,font=\scriptsize] {$2R$} (2.5,5.85);
    \node[pt,red!70!black,minimum size=4.4pt] at (2.5,2.5) {};
    \node at (2.5,-1.1) {(a)~Mitra--Virk};
  \end{scope}
  \begin{scope}[xshift=6.572cm]
    \begin{scope}
    \clip (-2,0) rectangle (4,6.3);
    \path[fill=gray!15] (-2,0) rectangle (4,1.5);
    \draw[rung,dashed] (-2,1.5) -- (4,1.5);
    \draw[rung,dashed] (-2,3) -- (4,3);
    \draw[rung,dashed] (-2,4.5) -- (4,4.5);
    \draw[rung,dashed] (-2,6) -- (4,6);
    \draw[tent,fill=blue!45,fill opacity=0.12,draw opacity=0.55] (-1.5,1.5) -- (0,3) -- (1.5,1.5) -- (0,0) -- cycle;
    \draw[tent,fill=blue!45,fill opacity=0.12,draw opacity=0.55] (0,1.5) -- (1.5,3) -- (3,1.5) -- (1.5,0) -- cycle;
    \draw[tent,fill=blue!45,fill opacity=0.12,draw opacity=0.55] (-0.75,3) -- (0.75,4.5) -- (2.25,3) -- (0.75,1.5) -- cycle;
    \node[pt,minimum size=2.6pt] at (-1.5,1.5) {};
    \node[pt,minimum size=2.6pt] at (0,1.5) {};
    \node[pt,minimum size=2.6pt] at (1.5,1.5) {};
    \node[pt,minimum size=2.6pt] at (3,1.5) {};
    \node[pt,minimum size=2.6pt] at (-0.75,3) {};
    \node[pt,minimum size=2.6pt] at (0.75,3) {};
    \node[pt,minimum size=2.6pt] at (2.25,3) {};
    \node[pt,minimum size=2.6pt] at (3.75,3) {};
    \node[pt,minimum size=2.6pt] at (-1.5,4.5) {};
    \node[pt,minimum size=2.6pt] at (0,4.5) {};
    \node[pt,minimum size=2.6pt] at (1.5,4.5) {};
    \node[pt,minimum size=2.6pt] at (3,4.5) {};
    \node[pt,minimum size=2.6pt] at (-0.75,6) {};
    \node[pt,minimum size=2.6pt] at (0.75,6) {};
    \node[pt,minimum size=2.6pt] at (2.25,6) {};
    \node[pt,minimum size=2.6pt] at (3.75,6) {};
    \end{scope}
    \draw[thick,-{Stealth[length=2mm]}] (-2,0) -- (4.4,0) node[right,inner sep=2pt] {$z$};
    \draw[ax] (-2,0) -- (-2,6.9) node[above,inner sep=1pt] {$\operatorname{pers}$};
    \draw (-2,1.5) -- (-2.2,1.5) node[left,inner sep=1.5pt,font=\scriptsize] {$\tau/2$};
    \draw (-2,3) -- (-2.2,3) node[left,inner sep=1.5pt,font=\scriptsize] {$\tau$};
    \draw (-2,4.5) -- (-2.2,4.5) node[left,inner sep=1.5pt,font=\scriptsize] {$3\tau/2$};
    \draw (-2,6) -- (-2.2,6) node[left,inner sep=1.5pt,font=\scriptsize] {$2\tau$};
    \draw[{Stealth[length=1.3mm]}-{Stealth[length=1.3mm]}] (-0.75,6.35) -- node[above,inner sep=1pt,font=\scriptsize] {$\tau/2$} (0.75,6.35);
    \draw[red!70!black,thin] (-0.25,1.2) rectangle (1.75,3.25) node[above right,inner sep=1pt,font=\scriptsize] {(c)};
    \node[pt,red!70!black,minimum size=4.4pt] at (1,1.9) {};
    \node at (1,-1.1) {(b)~Lemma~\ref{lem:config}};
  \end{scope}
  \begin{scope}
    \clip (17.6,0.35) rectangle (22.2,5.065);
    \draw[rung,dashed] (17.6,1.04) -- (22.2,1.04);
    \draw[rung,dashed] (17.6,4.49) -- (22.2,4.49);
    \draw[tent,fill=blue!45,fill opacity=0.12,draw opacity=0.55] (14.725,1.04) -- (18.175,4.49) -- (21.625,1.04) -- (18.175,-2.41) -- cycle;
    \draw[tent,fill=blue!45,fill opacity=0.12,draw opacity=0.55] (18.175,1.04) -- (21.625,4.49) -- (25.075,1.04) -- (21.625,-2.41) -- cycle;
    \draw[tent,fill=blue!45,fill opacity=0.12,draw opacity=0.55] (16.45,4.49) -- (19.9,7.94) -- (23.35,4.49) -- (19.9,1.04) -- cycle;
  \end{scope}
  \draw[red!70!black,thin] (17.6,0.35) rectangle (22.2,5.065);
  \node[pt,minimum size=3pt] at (18.175,1.04) {};
  \node[pt,minimum size=3pt] at (21.625,1.04) {};
  \node[pt,minimum size=3pt] at (19.9,4.49) {};
  \draw[red!70!black,densely dotted,thick] (20.475,1.96) -- (20.475,1.04) node[midway,left,inner sep=1.5pt,font=\tiny] {$\delta$};
  \draw[red!70!black,thick] (20.475,1.04) -- (21.625,1.04) node[midway,below,inner sep=1.5pt,font=\tiny] {$h$};
  \draw[red!70!black,densely dotted,thick] (20.475,1.96) -- (20.475,4.49) node[midway,right,inner sep=1.5pt,font=\tiny] {$\tau/2-\delta$};
  \draw[red!70!black,thick] (19.9,4.49) -- (20.475,4.49) node[midway,above,inner sep=1.5pt,font=\tiny] {$\tau/4-h$};
  \node[pt,red!70!black,minimum size=4.4pt] at (20.475,1.96) {};
  \node[right,font=\scriptsize] at (22.2,1.04) {$P_j$};
  \node[right,font=\scriptsize] at (22.2,4.49) {$P_{j+1}$};
  \node at (19.9,-1.1) {(c)~Zoom of (b)};
\end{tikzpicture}
\caption{The covers of Mitra--Virk ($R=1$) and of Lemma~\ref{lem:config} ($\tau=3$) in the coordinates $z=(b+d)/2$, $\operatorname{pers}=(d-b)/2$ of a bar $[b,d)$, in which $\ell^\infty$ is $|z-z'|+|\operatorname{pers}-\operatorname{pers}'|$: the diagonal is the $z$-axis, the tent supports, of radius $\tfrac32$, are diamonds, and the shaded strip supports the diagonal tent.
Centers (dots) lie in rows (dashed) with a common step between neighboring centers, alternate rows shifted by half a step.
(a)~Rows $R$ apart, half the step $2R$: four tents are positive at the marked point.
(b)~Rows $\tau/2$ apart, the full step: at most three are.
(c)~The multiplicity count in the proof: the marked point is offset by $\delta$ and $h$ from the nearest center of $P_j$, hence by $\tau/2-\delta$ and $\tau/4-h$ from that of $P_{j+1}$.}
\label{fig:covers}
\end{figure}

\begin{lemma}[Bounded-multiplicity covers of $\D{n}$]\label{lem:config}
For every $\tau>0$, $(\D{n},\db)$ carries a countable $(\tau,\tfrac14)$-admissible cover of multiplicity $\le \mathfrak m_n$, with $\mathfrak m_n$ depending only on $n$; one may take $\mathfrak m_n=3^n$.
\end{lemma}

\begin{proof}
\emph{Coordinates.}
Describe $x=(b,d)\in H$, a bar $[b,d)$, by its midpoint $z(x)=(b+d)/2$ and half-length $\operatorname{pers}x=(d-b)/2=\db(x,\Delta)$, so that $\db(x,y)=\min\{\norm{x-y}_\infty,\max(\operatorname{pers}x,\operatorname{pers}y)\}$ by Definition~\ref{def:config}.
Since $\max(|u-v|,|u+v|)=|u|+|v|$, the $\ell^\infty$-distance of two bars is the shift of their midpoints plus the change of their half-lengths,
\[
  \norm{x-y}_\infty=|z(x)-z(y)|+|\operatorname{pers}x-\operatorname{pers}y|,
\]
so in the coordinates $(z,\operatorname{pers})$ of Figure~\ref{fig:covers} the diagonal is the horizontal axis and $\ell^\infty$-balls are diamonds.

\emph{The cover, $n=1$.}
The \emph{row} $P_j$, $j\ge0$, consists of the bars of half-length $(j+1)\tau/2$ whose midpoints lie in $\tfrac\tau2\Z$, shifted by $\tfrac\tau4$ when $j$ is odd.
The cover consists of the \emph{diagonal tent} $\varphi_{\Delta,\tau/2}$ and the tents $\varphi_{p,\tau/2}$, $p\in P:=\bigcup_jP_j$, a countable family.
Since $\operatorname{pers}p\ge\tau/2$, the cap in $\db(p,\cdot)$ never binds below $\tau/2$: $\varphi_{p,\tau/2}(y)=\max\{\tau/2-\norm{p-y}_\infty,0\}$ on $H$, and $\varphi_{p,\tau/2}(\Delta)=0$.
For a bar $x$ and a row, let $h\in[0,\tau/4]$ be the horizontal distance $|z(x)-z(p)|$ to the nearest center $p$ of the row; consecutive rows are offset by half their step $\tau/2$, so in the next row the nearest center is at horizontal distance $\tau/4-h$.

\emph{Admissibility: some tent is at least $\tau/8$ at every $x$.}
If $\operatorname{pers}x\le3\tau/8$, the diagonal tent is $\tau/2-\operatorname{pers}x\ge\tau/8$.
If $\operatorname{pers}x\in[3\tau/8,\tau/2]$, the nearest center of $P_0$ is at distance $h+\tau/2-\operatorname{pers}x\le3\tau/8$.
Otherwise $x$ lies between $P_j$ and $P_{j+1}$, at vertical distances $\delta$ and $\tau/2-\delta$, and the nearest centers of the two rows are at distances $h+\delta$ and $(\tau/4-h)+(\tau/2-\delta)$, which add up to $3\tau/4$; one is at most $3\tau/8$.
So $\lambda_0\ge\tau/8$, and the cover is $(\tau,\tfrac14)$-admissible.

\emph{Multiplicity: at most $3$ tents are positive at every $y$.}
A row at vertical distance $\delta$ from $y$ contributes its centers at horizontal distance less than $\tau/2-\delta$, among $h,\ \tau/2-h,\ \tau/2+h,\dots$; so only rows with $\delta<\tau/2$ contribute, at most $2$ each.
If $\operatorname{pers}y<\tau/2$, only $P_0$ contributes, and the diagonal tent makes $3$.
If $\operatorname{pers}y\ge\tau/2$, the diagonal tent vanishes, and only the bracketing rows contribute: $P_j$, at distance $\delta\in[0,\tau/2)$, gives $[h<\tau/2-\delta]+[h>\delta]$, and $P_{j+1}$, at distance $\tau/2-\delta$, gives $[\tau/4-h<\delta]+[\tau/4+h<\delta]$, where $[\,\cdot\,]$ is $1$ if the condition holds and $0$ otherwise.
Four would need $h>\delta>\tau/4+h$.

\emph{General $n$.}
The centers are the $n$-point multisets $T$ drawn from $P\cup\{\Delta\}$, with radius $\tau/2$.
Given $A=\{a_1,\dots,a_n\}$, the case $n=1$ gives $t_i$ with $\db(t_i,a_i)\le3\tau/8$, and the matching $a_i\mapsto t_i$ gives $\varphi_{T,\tau/2}(A)\ge\tau/8$ for $T=\{t_i\}$.
If $\varphi_{T,\tau/2}(A)>0$, some matching pairs each $a_i$ with a point $t_i$ of $T$ with $\varphi_{t_i,\tau/2}(a_i)>0$; each $a_i$ has at most $3$ such $t_i$, so at most $3^n=:\mathfrak m_n$ multisets $T$ qualify.
\end{proof}

\noindent The covers of Lemma~\ref{lem:config} exist at every scale, so Proposition~\ref{prop:covers-asdim} gives $\asdim\D n\le3^n-1$.
For $n=1$ this is the exact value $2$, and multiplicity $3$ is the least possible for a multi-scale admissible family on $\D1$: far from the diagonal $\D1$ contains isometric copies of arbitrarily large $\ell^\infty$-squares, so $\asdim\D1\ge2$ by Lemma~\ref{lem:cube-asdim}, while a family of multiplicity $\mathfrak m$ gives $\asdim\D1\le\mathfrak m-1$ by Proposition~\ref{prop:covers-asdim}.
For $n\ge2$ the bound is far from Mitra and Virk's exact value $2n$~\cite{Mitra2021,MitraVirk2024corrigendum}, because the multiplicities of the $n$-tuple covers multiply.
The ratio $\tfrac14$ is within a factor $2$ of the best possible, since the diagonal caps it: no tent cover of $\D1$ of finite multiplicity is $(\tau,\theta)$-admissible with $\theta>\tfrac12$.
Indeed, by finite multiplicity only finitely many tents are positive at $\Delta$; call their centers $c_1,\dots,c_k$, and take the point $x=(s,s+\tau/2)$, of persistence $\tau/4$, with $s$ so large that $s\ge b_i+\tau/4$ for every $c_i=(b_i,d_i)\ne\Delta$, so that $\norm{x-c_i}_\infty\ge\tau/4$.
A tent vanishing at $\Delta$ is at most $\operatorname{pers}x=\tau/4$ at $x$, since it is $1$-Lipschitz and $x$ is at distance $\operatorname{pers}x$ from $\Delta$.
A tent positive at $\Delta$, with center $c_i$ and radius $r\le\tau/2$, is at most $r-\operatorname{pers}x\le\tau/4$ at $x$, since $\db(c_i,x)\ge\operatorname{pers}x$: for $c_i=\Delta$ with equality, and otherwise by Definition~\ref{def:config}, as $\norm{x-c_i}_\infty\ge\operatorname{pers}x$.
So every tent is at most $\tau/4$ at $x$, whence $\lambda_0\le\tau/4$ and $\theta\le\tfrac12$.

\begin{theorem}[Coarse embedding of $n$-point diagrams]
\label{thm:diagrams}
For every $\tau>0$ there is a $1$-Lipschitz map $\Phi\colon(\D{n},\db)\to\ell^2$, built from tents centered at $n$-point diagrams, with
\[
  \norm{\Phi(A)-\Phi(B)}\ \ge\ \tfrac1{16}\,3^{-n/2}\sqrt{\tau\,\db(A,B)}\qquad\text{whenever }\db(A,B)\ge\tau.
\]
In particular $\Phi$ is a coarse embedding.
\end{theorem}

\begin{proof}
The covers of Lemma~\ref{lem:config} at the scales $s_j=\tau2^j$ form a multi-scale admissible family $\mathcal F$ on $(\D{n},\db)$ of ratio $\tfrac14$ and multiplicity $\le3^n$, whose tents $\varphi_{T,r}(A)=\max\{r-\db(T,A),0\}$ are centered at $n$-point diagrams $T$, the \emph{$n$-tuple landmarks}.
With the dyadic weights, Theorem~\ref{thm:multi}\,(i) makes $\Phi=\Phi(\cdot;\mathcal F)$ $1$-Lipschitz with lower gauge $c_{1/4}\,3^{-n/2}\sqrt{\tau t}=\tfrac1{16}\,3^{-n/2}\sqrt{\tau t}$ for $t\ge\tau$.
The same argument applies to every multi-scale admissible family on $(\D n,\db)$ of ratio $\theta$ and multiplicity $\le\mathfrak m$, with lower gauge $c_\theta\,\mathfrak m^{-1/2}\sqrt{\tau t}$; its covers need not be lattices.
\end{proof}

\begin{remark}[Relation to Mitra--Virk]\label{rem:mv}
Mitra and Virk~\cite{Mitra2024} embed $(\D{n},\db)$ by the same construction: tents centered at the $n$-point diagrams of a fixed lattice, stacked over scales $R_k$ with weights $w_k$.
Their per-scale floor~\cite[Lem.~3.18]{Mitra2024} corresponds to the two-tent witness in the proof of Theorem~\ref{thm:multi}, and their coarse and uniform embeddings~\cite[Thm.~4.3]{Mitra2024} are the same construction on other schedules (Remark~\ref{rem:schedules}).
The constructions differ in three respects.

\emph{(i)~The cover.} At scale $R$ their cover has tents of radius $3R/2$, height at least $R/8$ everywhere and multiplicity at most $4^n$~\cite[Def.~3.13, Lemmas~3.12 and~3.14]{Mitra2024}, so it is $(3R,\tfrac1{12})$-admissible; the cover of Lemma~\ref{lem:config} has ratio $\tfrac14$ and multiplicity $3^n$.
Both place their centers in rows parallel to the diagonal, alternate rows shifted by half the step between neighboring centers; their rows are half a step apart, forming a square grid, and ours a full step apart, which lowers the multiplicity from $4$ to $3$ (Figure~\ref{fig:covers}).
Whatever the schedule, the lower gauge of Theorem~\ref{thm:multi} carries the factor $\theta/\sqrt{\mathfrak m}$, here $\tfrac14\,3^{-n/2}$ against $\tfrac1{12}\,2^{-n}$, a factor $3(2/\sqrt3)^n$ in favor of Lemma~\ref{lem:config}.

\emph{(ii)~The gauge.} Both maps are $1$-Lipschitz once the weights of~\cite[Thm.~4.3]{Mitra2024} form a unit vector, as that theorem requires (their Example~4.6 takes $w_k=1/k$, which normalizes to $w_k=\sqrt6/(\pi k)$), and both start at the same base radius: in the proof of their Theorem~4.3 the scale-$R_k$ block has tents of radius $R_k/2$, the radius $3R/2$ of~(i) at $R=R_k/3$.
The gauge stated there is the step function $\sum_kc_nw_kR_k\chi_{[R_k,R_{k+1})}$, $c_n=2^{-n-3}/3$, at most $c_n\sqrt t$ on their schedule $R_k=k^2$ and $0$ below $R_1=1$, like the dyadic map of Theorem~\ref{thm:diagrams} below $\tau$; for $\tau=1$ the latter exceeds it on $[1,\infty)$ by at least the factor $\tfrac32(2/\sqrt3)^n$.
At small scales they add a second map, on scales tending to $0$, and combine the two; the two-sided ladder of Remark~\ref{rem:schedules} instead gives a single map whose gauge is positive at every $t>0$.
Summing the finer scales, as in~\cite[Rem.~4.5]{Mitra2024}, gives their map a gauge of order $c_nt^{3/4}$, which eventually exceeds the $\sqrt{\tau t}$ of the dyadic weights; the nearly linear gauge of Theorem~\ref{thm:multi}\,(ii) eventually exceeds both.

\emph{(iii)~Flexibility.} The proof of Theorem~\ref{thm:diagrams} admits any admissible family, so the landmarks may be placed where the data concentrate, for instance by farthest-point sampling~\cite{Gonzalez1985}, as in PALACE~\cite{PaperII}.
For the $p$-Wasserstein distance $W_p$ with the same ground metric and diagonal cost as $\db$ (Definition~\ref{def:config}), $\db\le W_p\le(2n)^{1/p}\db$ on $\D n$, a matching having at most $2n$ pairs (cf.~\cite[Prop.~3.1]{Mitra2021}); so each such map is also a $1$-Lipschitz coarse embedding of $(\D n,W_p)$, with lower gauge $\rho_-\bigl((2n)^{-1/p}t\bigr)$ for $t\ge(2n)^{1/p}\tau$.
\end{remark}

\subsection{Slicing to multiparameter modules}\label{sec:multiparam}

Multiparameter modules reduce to the diagram case one slice at a time.
Fix a \emph{slice net}, a finite family $\mathcal L_S=\{L_1,\dots,L_S\}\subset\mathcal L$ of normalized lines (Definition~\ref{def:dmatch}), and the sliced distance $\dM^{S}(M,N)=\max_k\db\bigl(\mathcal B(M\circ L_k),\mathcal B(N\circ L_k)\bigr)\le\dM(M,N)$; the normalization builds the weights $w(u)$ into the parametrization, so that no explicit weight appears.

Let $\RN$ be the class of finitely presented $d$-parameter modules with at most $N$ bars on every slice.
Restricting a finite presentation to a slice gives a finitely presented one-parameter module whose births and deaths are entry times~\eqref{eq:entry} of generators and relations.
Its barcode has bars $[b,d)$ with $d\in(b,\infty]$, the \emph{essential} bars $[b,\infty)$ included, so barcodes live not in $\D N$ but in
\[
  \D{\le N}:=\bigsqcup_{a=0}^{N}E_a\times\D{N-a},
\]
where $E_a$ is the set of $a$-point multisets in $\R$ under the matching distance $d_E$, with $d_E(A,B)=\min_\sigma\max_i|a_i-b_{\sigma(i)}|$, each piece carries the max metric, and distinct pieces are at infinite distance.
A barcode with $a$ essential bars and at most $N-a$ finite ones is the pair (essential births, finite bars), a point of $E_a\times\D{N-a}$; this identifies the barcodes with at most $N$ bars, under $\db$, isometrically with $\D{\le N}$.
Indeed, a matching of finite cost pairs essential bars with essential bars, at cost $|b-b'|$, and finite bars with finite bars or the diagonal; the two choices are independent, so $\db$ is the maximum of $d_E$ and of the bottleneck distance of the finite parts.
Sorting identifies $E_a$ isometrically with the cone $\{x_1\le\dots\le x_a\}\subset(\R^a,\norm{\cdot}_\infty)$.

\begin{lemma}[Bounded-multiplicity covers of $\D{\le N}$]\label{lem:config-ess}
For every $\tau>0$, $\D{\le N}$ carries a countable $(\tau,\tfrac14)$-admissible cover of multiplicity at most $3^N$.
\end{lemma}

\begin{proof}
On $E_a$, in sorted coordinates, take as centers the points of the lattice $\tfrac{3\tau}4\Z^a$ that lie in the cone, with radius $\tau/2$.
Rounding each coordinate of $x\in E_a$ to a nearest multiple of $3\tau/4$, ties downward, keeps the order and moves each coordinate by at most $3\tau/8$, so some tent is at least $\tau/8$ at $x$.
A tent is positive at $x$ only if each coordinate of its center is within $\tau/2$ of that of $x$, and an open interval of length $\tau$ contains at most two multiples of $3\tau/4$, so at most $2^a$ tents are positive at $x$.
On $\D{N-a}$ take the cover of Lemma~\ref{lem:config}, with radii $\tau/2$, $\lambda_0\ge\tau/8$ and multiplicity at most $3^{N-a}$ (a single tent when $a=N$, as $\D0$ is a point).
For centers $p\in E_a$ and $q\in\D{N-a}$ the tent of radius $\tau/2$ at $(p,q)$ in the max metric is $\varphi_{(p,q),\tau/2}(x,y)=\min\{\varphi_{p,\tau/2}(x),\varphi_{q,\tau/2}(y)\}$, so the product cover has $\lambda_0\ge\tau/8$ and multiplicity at most $2^a3^{N-a}\le3^N$.
A tent centered on one piece vanishes on the others, so the union over $a$ is $(\tau,\tfrac14)$-admissible of multiplicity at most $3^N$.
\end{proof}

Each restriction $M\circ L_k$ is a point of $(\D{\le N},\db)$; embedding $\D{\le N}$ piece by piece and stacking over the slices embeds the module.

\begin{theorem}[Coarse embedding of multiparameter modules]
\label{thm:modules}
For every slice net $\mathcal L_S$ and every $\tau>0$ there is a map $\Psi\colon\RN\to\ell^2$ such that, with $t=\dM^{S}(M,M')$,
\[
  c_N\sqrt{\tau t}\ \le\ \norm{\Psi(M)-\Psi(M')}\ \le\ \sqrt S\,t,
\]
the upper bound whenever $t<\infty$ and the lower bound when $t\ge\tau$, where $c_N=\tfrac1{16}\,3^{-N/2}$.
So $\Psi$ is a coarse embedding of each class of modules at finite sliced distance, with no bound on the parameter range, and the quotient of $(\RN,\dM^S)$ by $\dM^S=0$ has finite asymptotic dimension.
\end{theorem}

\begin{proof}
Let $\psi\colon(\D{\le N},\db)\to\ell^2$ be, on each piece $E_a\times\D{N-a}$ and with a base point of its own, the embedding of Theorem~\ref{thm:multi}\,(i) built from the covers of Lemma~\ref{lem:config-ess} at the scales $s_j$.
On each piece $\psi$ is $1$-Lipschitz with lower gauge $\rho_-(t)=c_N\sqrt{\tau t}$ for $t\ge\tau$, since $c_{1/4}=\tfrac1{16}$ and the multiplicity is at most $3^N$.
Put $\Psi(M)=(\psi(M\circ L_k))_{k=1}^S\in(\ell^2)^S\cong\ell^2$.

Modules at finite sliced distance are those with the same number of essential bars on every slice of the net; this partitions $\RN$ into at most $(N+1)^S$ classes; for $M,M'$ in one class, $M\circ L_k$ and $M'\circ L_k$ lie in one piece of $\D{\le N}$ for every $k$.
With $g_k=\norm{\psi(M\circ L_k)-\psi(M'\circ L_k)}$ and $\delta_k=\db(M\circ L_k,M'\circ L_k)$ we have $g_k\le\delta_k$, $g_k\ge\rho_-(\delta_k)$ when $\delta_k\ge\tau$, and $\norm{\Psi(M)-\Psi(M')}^2=\sum_k g_k^2$.
If $k^\star$ attains $\delta_{k^\star}=\max_k\delta_k=t$, the block $k^\star$ gives the lower bound, while $\sum_k g_k^2\le\sum_k\delta_k^2\le St^2$ gives the upper.
Since $\rho_-\to\infty$, $\Psi$ is a coarse embedding on each class; it sees $M$ only through its fibered barcode over $\mathcal L_S$, so its lower gauge is in terms of $\dM^S$, which is at most $\dM$.
For the last assertion, the slice map $M\mapsto(\mathcal B(M\circ L_k))_{k\le S}$ embeds the quotient isometrically into the max-product of $S$ copies of $(\D{\le N},\db)$, which has finite asymptotic dimension by Lemma~\ref{lem:config-ess}, Proposition~\ref{prop:covers-asdim} and the product theorem~\cite{BellDranishnikov2008}; sharp bounds are in the companion paper~\cite{MajhiMitraVirk-asdim}.
\end{proof}

\begin{remark}[Computation]\label{rem:complex}
\emph{Whole-tuple landmarks.} The embedding scores a diagram $A$ against each landmark $T$ by one bottleneck matching, $\db(T,A)$, and its lower gauge holds for every pair at distance at least $\tau$.
The embedding COMPLEX of \cite{Majhi2026-complex} scores a diagram by a sum of one-point tents instead, which is cheaper but many-to-one, so its lower gauge rests on a \emph{coherence} hypothesis checked on the data.

\emph{Bounded frames.} For modules whose generators and relations lie in a box of side $L$, every normalized line sees all births and deaths in one interval of length at most $L$, because the entry time~\eqref{eq:entry} along a normalized line is $1$-Lipschitz in the grade.
Pairs at finite distance are then at most $L$ apart, so the scales $s_j\le L$, about $\log_2(L/\tau)$ of them, suffice: truncating the stack there keeps the lower gauge and only decreases the upper bound.
The covers of Lemmas~\ref{lem:config} and~\ref{lem:config-ess} have their centers on lattices, so only finitely many tents of each scale meet the frame, and the stack lands in a finite-dimensional Euclidean space, the bounded-domain regime of Mitra--Virk~\cite[Thm.~5.1]{Mitra2024}.
\end{remark}

\section{Bounded presentations and slices}\label{sec:slices}

The results of Section~\ref{sec:persistence} are stated for modules with at most $N$ bars per slice, a condition inherited from one parameter.
This section compares it with the size of a presentation: the number of generators bounds the number of bars per slice (Lemma~\ref{lem:slice-cardinality}), but a bound on bars per slice does not control the coarse geometry, since a single bar per slice is compatible with infinite asymptotic dimension (Theorem~\ref{thm:staircase}).
Throughout this section modules are $\R^2$-indexed.

Presentations, minimal resolutions and the Betti positions $\beta_0^M,\beta_1^M$ are as in Section~\ref{sec:modules}.

\begin{definition}\label{def:Mnm}
For $n,m\in\mathbb{N}$ let $\Mnm$ denote the set of isomorphism classes
$[M]$ of finitely presented $\mathbb{R}^2$-indexed persistence modules
whose minimal free presentation $\partial_1\colon F_1\to F_0$ satisfies
\[
  \operatorname{rank}F_0\le n
  \qquad\text{and}\qquad
  \operatorname{rank}F_1\le m;
\]
equivalently, $|\beta_0^M|\le n$ and $|\beta_1^M|\le m$, counted with
multiplicity.  The condition is well posed because the minimal
presentation is unique up to isomorphism.
The class $\RN$ of Theorem~\ref{thm:modules}, the finitely presented modules with at most $N$ bars on every slice (Section~\ref{sec:multiparam}), is now taken with $d=2$.
\end{definition}

\subsection{Slice cardinality}

Restricted to a slice, a finitely presented module is a p.f.d.\ one-parameter module with a barcode (Section~\ref{sec:modules}), and the number of bars is controlled by the number of generators alone.

\begin{lemma}[Slice cardinality]\label{lem:slice-cardinality}
Let $M\in\Mnm$ and let $L$ be a line of positive slope, parametrized at unit speed or as in $\mathcal L$.  Then
$|\mathcal{B}(M\circ L)|\le n$, and every bar has the form $[b,d)$ with $d\in(b,\infty]$.
\end{lemma}

\begin{proof}
Restriction along $L$ is exact, being pointwise.
A free module $k(-\gamma)$ restricts to the interval module on the closed ray $\{t:\gamma\le L(t)\}=[e_L(\gamma),\infty)$, with the entry time $e_L(\gamma)=\max_i(\gamma_i-b_i)/v_i$ as in~\eqref{eq:entry}.
Hence $M\circ L=\operatorname{coker}(F_1\circ L\to F_0\circ L)$ is a finitely presented one-parameter module, and so a finite sum of intervals $[b_j,d_j)$ with $d_j\in(b_j,\infty]$.
For a one-parameter module $V$ put $V_{<t}:=\sum_{s<t}\operatorname{im}(V(s)\to V(t))$; each bar $[b,d)$ contributes exactly $1$ to $\dim V(t)/V_{<t}$ at $t=b$ and $0$ elsewhere, so the number of bars is $\sum_t\dim V(t)/V_{<t}$.
The surjection $F_0\circ L\to M\circ L$ maps $(F_0\circ L)_{<t}$ into $(M\circ L)_{<t}$ and hence induces surjections $(F_0\circ L)(t)/(F_0\circ L)_{<t}\to(M\circ L)(t)/(M\circ L)_{<t}$.
The left side has dimension $\#\{i:e_L(\gamma_i)=t\}$; summing over $t$ gives at most $\operatorname{rank}F_0\le n$ bars.
\end{proof}

\noindent Relations only shorten or delete bars, so the bound is independent of $m$: $\Mnm\subseteq\mathcal R_n$ for every $m$, an inclusion far from equality by Theorem~\ref{thm:staircase}.

\begin{remark}[Other distances]
Interleavings with independent shifts in the two parameters, the generalized interleavings of Bubenik, de Silva and Scott~\cite{BubenikDeSilvaScott2015}, give distances each bi-Lipschitz equivalent to its counterpart, $\di$ or $\dM$, so nothing here depends on that choice; the $\ell^p$-distances of Bjerkevik and Lesnick~\cite{BjerkevikLesnick2021} are of Wasserstein type and, for $p<\infty$, not bi-Lipschitz to $\di$.
\end{remark}

\subsection{Staircases: bounding bars per slice is not enough}\label{sec:staircase}

The hypothesis of Section~\ref{sec:persistence}, at most $N$ bars on every slice, is the one-parameter condition carried over verbatim.
It does not control the presentation, and it does not control the coarse geometry.
Write $\one=(1,1)$.

\begin{definition}[Staircase modules]\label{def:staircase}
Let $A=\{u_1,\dots,u_q\}\subset\R^2$ be an antichain, indexed so that $u_1^1<\dots<u_q^1$ and $u_1^2>\dots>u_q^2$, let $U_A=\bigcup_i(u_i+\R_{\ge0}^2)$ be the closed upset it generates, and let $k[U_A]$ be its indicator module: $k$ on $U_A$, zero elsewhere, identity structure maps on $U_A$.
\end{definition}

\begin{lemma}\label{lem:staircase-presentation}
$k[U_A]$ is finitely presented with generators at $u_1,\dots,u_q$ and relations at the $q-1$ joins $u_i\vee u_{i+1}$; so $k[U_A]\in\mathcal M_2^{q,q-1}$.
Every slice meets $U_A$ in a closed ray, so $k[U_A]$ has exactly one bar, essential, on every slice: $k[U_A]\in\mathcal R_1$ for every $q$.
\end{lemma}

\begin{proof}
Write $I(x)=\{i:u_i\le x\}$ for $x\in\R^2$ and $u_{i,i+1}=u_i\vee u_{i+1}$.
The set $I(x)$ is an interval of indices: if $u_i\le x$ and $u_j\le x$ with $i<l<j$, then $u_l^1\le u_j^1\le x^1$ and $u_l^2\le u_i^2\le x^2$.
Let $F_0=\bigoplus_i k(-u_i)$, $F_1=\bigoplus_{i<q}k(-u_{i,i+1})$ and $\partial(e_{i,i+1})=e_i-e_{i+1}$.
At $x$ the cokernel is $k^{I(x)}$ modulo the differences of consecutive present indices, since $u_{i,i+1}\le x$ exactly when both are present; this is $k$ if $I(x)$ is nonempty and $0$ otherwise, with identity structure maps, which is $k[U_A]$.
It remains to check that the presentation is minimal, that is, $\beta_0^M=\{u_1,\dots,u_q\}$ and $\beta_1^M=\{u_{i,i+1}:i<q\}$ for $M=k[U_A]$.
A grade $x$ carries a generator of a module $N$ exactly when $\sum_{y<x}\operatorname{im}\bigl(N(y)\to N(x)\bigr)\neq N(x)$, and $\beta_1^M=\beta_0^Z$ for $Z=\ker(F_0\to M)$ once $F_0$ is minimal (Section~\ref{sec:modules}).
For $\beta_0^M$: if $y<u_i$ then $y\notin U_A$, since $u_j\le y<u_i$ would contradict $A$ being an antichain; so $M(y)=0$ for all $y<u_i$ and $u_i$ carries a generator.
At any other $x\in U_A$ there are $i$ and $y$ with $u_i\le y<x$, and $M(y)\to M(x)$ is the identity of $k$; so $x$ carries none.
For $\beta_1^M$: at $x$ the map $F_0(x)=k^{I(x)}\to M(x)$ is the sum of coordinates when $I(x)\ne\emptyset$, so $Z(x)=\{v\in k^{I(x)}:\sum_i v_i=0\}$ is spanned by the $e_i-e_{i+1}$ with $i,i+1\in I(x)$.
Now $I(u_{i,i+1})=\{i,i+1\}$: for $l<i$, $u_l\le u_{i,i+1}$ would force $u_l^2\le u_i^2$, and for $l>i+1$ it would force $u_l^1\le u_{i+1}^1$, both against the indexing.
Hence $Z(u_{i,i+1})=k\,(e_i-e_{i+1})$, while $|I(y)|\le1$ and $Z(y)=0$ for every $y<u_{i,i+1}$, because $u_i\le y$ and $u_{i+1}\le y$ together give $y\ge u_{i,i+1}$; so $u_{i,i+1}$ carries a generator of $Z$.
At any other $x$ with $Z(x)\neq0$ each spanning vector $e_i-e_{i+1}$ has $u_{i,i+1}\le x$, $u_{i,i+1}\ne x$, so lies in the image of $Z(u_{i,i+1})$ with $u_{i,i+1}<x$; so $x$ carries none.
Thus $\beta_0^Z=\{u_{1,2},\dots,u_{q-1,q}\}$ and $\partial$ is a minimal presentation.
A line of positive slope enters the closed upset $U_A$ once and never leaves it.
\end{proof}

\begin{lemma}[Interleavings between upset modules]\label{lem:upset-interleaving}
For nonempty closed upsets $U,V\subset\R^2$,
\[
  \di(k[U],k[V])=\inf\{\eps\ge0:\ U+\eps\one\subseteq V\ \text{and}\ V+\eps\one\subseteq U\}.
\]
If $U=U_A$ and $V=U_B$ with $A=\{u_i\}$, $B=\{v_j\}$, then $U+\eps\one\subseteq V$ iff for every $i$ some $j$ has $v_j\le u_i+\eps\one$, so
\[
  \di(k[U_A],k[U_B])=\max\Bigl\{\max_i\min_j\|(v_j-u_i)_+\|_\infty,\ \max_j\min_i\|(u_i-v_j)_+\|_\infty\Bigr\},
\]
with $(\cdot)_+$ the componentwise positive part.
Moreover $\dM(k[U],k[V])=\di(k[U],k[V])$.
\end{lemma}

\begin{proof}
Since $k[V][\eps]=k[V-\eps\one]$, a morphism $f\colon k[U]\to k[V][\eps]$ is a family of scalars $f_x$, $x\in U\cap(V-\eps\one)$, zero elsewhere, with $f_y=f_x$ whenever $x\le y$ in $U\cap(V-\eps\one)$ and $f_y=0$ whenever $y\ge x$ for some $x\in U\setminus(V-\eps\one)$.
If $(f,g)$ is an $\eps$-interleaving, the composite $g[\eps]\circ f\colon k[U]\to k[U-2\eps\one]$ is the shift map, which is the identity of $k$ at every $x\in U$ (as $U\subseteq U-2\eps\one$); so $f_x\ne0$ for every $x\in U$, forcing $U\subseteq V-\eps\one$, that is, $U+\eps\one\subseteq V$, and symmetrically $V+\eps\one\subseteq U$.
Conversely, if both inclusions hold, $f_x=1$ on $U$ and $g_y=1$ on $V$ are morphisms and their composites are the shifts.
For the corner formula, $u_i+\eps\one\in U_B$ iff $v_j\le u_i+\eps\one$ for some $j$ iff $\min_j\max_m(v_j^m-u_i^m)\le\eps$, and $U_A+\eps\one\subseteq U_B$ iff this holds for every $i$, since $U_A+\eps\one$ is generated by the $u_i+\eps\one$.

For the matching distance, $\dM\le\di$ (Definition~\ref{def:dmatch}).
For the reverse inequality use the slope-one lines $L_b(t)=b+t\one\in\mathcal L$: the slice $k[U]\circ L_b$ is the single bar $[e_U(b),\infty)$ with entry time $e_U(b)=\inf\{t:b+t\one\in U\}$, so $\db(\mathcal B(k[U]\circ L_b),\mathcal B(k[V]\circ L_b))=|e_U(b)-e_V(b)|$.
Suppose $\sup_b|e_U(b)-e_V(b)|\le\eta$.
For $x\in U$ take $b=x$: then $e_U(x)\le0$, so $e_V(x)\le\eta$ and $x+\eta\one\in V$ as $V$ is a closed upset; hence $U+\eta\one\subseteq V$, and symmetrically, so $\di(k[U],k[V])\le\eta\le\dM(k[U],k[V])$.
\end{proof}

\begin{lemma}[Cubes force dimension]\label{lem:cube-asdim}
Let $K\ge1$.
If a metric space $X$ contains, for every $L>0$, an isometric copy of the cube $([0,L]^K,\ell^\infty)$, then $\asdim X\ge K$ (cf.~\cite[Lemma~2.10]{Mitra2021}).
\end{lemma}

\begin{proof}
Suppose $\asdim X\le K-1$.
By the coloring characterization of asymptotic dimension~\cite{BellDranishnikov2008} there are families $\mathcal U_0,\dots,\mathcal U_{K-1}$ of subsets of $X$, of diameter at most some $D$, whose union covers $X$, such that distinct members of one family are at distance at least $3$.
Take $L>D+2$ and an isometric copy $Q$ of $([0,L]^K,\ell^\infty)$ in $X$.
For $U\in\mathcal U_j$ meeting $Q$ put $U'=\{x\in Q:d(x,U)\le1\}$: a closed subset of $Q$ of diameter at most $D+2<L$, so $U'$ meets no two opposite faces of the cube.
Distinct members of one family give sets $U'$ at distance at least $1$, hence pairwise disjoint; and since $Q$ is compact, only finitely many members of each family meet $Q$ (a point chosen in each would form an infinite $1$-separated subset of $Q$).
So the sets $U'$ form a finite closed cover of $Q$ by sets none of which meets two opposite faces, and by Lebesgue's covering theorem~\cite[Thm.~IV~2]{HurewiczWallman1941} some point of $Q$ lies in at least $K+1$ of them.
But a point lies in at most one $U'$ per family, hence in at most $K$ of them, a contradiction.
\end{proof}

\begin{theorem}[Staircases]\label{thm:staircase}
For $q\ge1$ and $R>0$ let $\mathcal S_q(R)$ be the family of staircase modules $k[U_A]$ with corners
\[
  u_i=\bigl(iR,(q-i)R\bigr)+(s_i,t_i),\qquad (s_i,t_i)\in[0,R/3]^2,\quad i=1,\dots,q .
\]
Then $(s_i,t_i)_{i}\mapsto k[U_A]$ is an isometry from the $\ell^\infty$-cube $([0,R/3]^2)^q\subset\R^{2q}$ onto $(\mathcal S_q(R),\di)=(\mathcal S_q(R),\dM)$.
Consequently, for every $N\ge1$,
\[
  \asdim(\RN,\di)=\asdim(\RN,\dM)=\infty,
\]
and neither $(\RN,\di)$ nor $(\RN,\dM)$ is doubling.
\end{theorem}

\begin{proof}
Since $s_i,t_i\in[0,R/3]$ and $R/3<R$, the corners satisfy $u_1^1<\dots<u_q^1$ and $u_1^2>\dots>u_q^2$, so each $A$ is an antichain indexed as in Definition~\ref{def:staircase} and Lemma~\ref{lem:staircase-presentation} applies.
Let $A=\{u_i\}$, $B=\{v_i\}$ be two corner configurations of the family.
For $j>i$ the first coordinate gives $v_j^1-u_i^1\ge(j-i)R-R/3\ge2R/3$, and for $j<i$ the second gives $v_j^2-u_i^2\ge2R/3$, while $\|(v_i-u_i)_+\|_\infty\le\|v_i-u_i\|_\infty\le R/3$.
So in Lemma~\ref{lem:upset-interleaving} the inner minimum is attained at $j=i$, and
\[
  \di(k[U_A],k[U_B])=\max_i\max\bigl\{\|(v_i-u_i)_+\|_\infty,\|(u_i-v_i)_+\|_\infty\bigr\}=\max_i\|u_i-v_i\|_\infty,
\]
the $\ell^\infty$-distance of the configurations in $\R^{2q}$; the same holds for $\dM$ by the lemma.
By Lemma~\ref{lem:staircase-presentation} the family lies in $\mathcal R_1\subseteq\RN$ and in $\mathcal M_2^{q,q-1}$.

Thus $(\RN,\di)$ and $(\RN,\dM)$ contain, for every $R$, an isometric copy of $([0,R/3]^{2q},\ell^\infty)$, and Lemma~\ref{lem:cube-asdim} gives $\asdim\ge2q$.
Hence $\asdim(\RN,\di)\ge2q$ for every $q$, and likewise for $\dM$.
For the failure of doubling, let $c$ be the center of the cube and fix $r\in(R/6,R/3)$.
The open ball $B(c,r)$ of $(\RN,\di)$ contains the $2^{2q}$ vertices of the cube, which lie at distance $R/6$ from $c$ and pairwise at distance $R/3$; a ball of radius $r/2$ has diameter at most $r<R/3$ and so contains at most one vertex, hence covering $B(c,r)$ takes at least $2^{2q}$ balls of radius $r/2$.
A doubling constant would bound this independently of $q$.
The same holds for $\dM$, the two distances agreeing on the family.
\end{proof}

\begin{corollary}[No coarse embedding]\label{cor:staircase-nonembed}
The space $\mathcal R_1$ contains every finite metric space isometrically, under $\di$ and under $\dM$.
Consequently, for every $N\ge1$, neither $(\RN,\di)$ nor $(\RN,\dM)$ coarsely embeds into Hilbert space.
\end{corollary}

\begin{proof}
A finite metric space $(F,d_F)$ with $p$ points embeds isometrically into $(\R^p,\ell^\infty)$ by the Fr\'echet map $x\mapsto(d_F(x,y))_{y\in F}$, whose image lies in the cube $[0,\diam F]^p$.
Padding with constant coordinates embeds it isometrically into the cube $([0,R/3]^2)^q$ whenever $2q\ge p$ and $R\ge3\diam F$, hence into $\mathcal S_q(R)\subseteq\mathcal R_1$ by Theorem~\ref{thm:staircase}.
A metric space $X$ that contains every finite metric space isometrically admits no coarse embedding into Hilbert space~\cite[Cor.~4.7(1)]{Mitra2021}: such an embedding would restrict, on isometric copies in $X$ of the finite metric spaces $(\mathbb Z_n)^m$ of Theorem~\ref{thm:dgly}, to maps with controls independent of $n$ and $m$, and the proof of Theorem~\ref{thm:dgly} shows that no such maps exist.
\end{proof}

The device, every finite metric space placed isometrically and far apart, is that of Bubenik and Wagner~\cite{BubenikWagner2020} and, independently, of Mitra and Virk~\cite[Lemma~4.1, Cor.~4.2 and Thm.~4.3]{Mitra2021} for the space of finite diagrams; here a single bar per slice suffices.
The same mechanism gives Zava's non-embeddability of the Gromov--Hausdorff space of finite metric spaces~\cite{Zava2025}, following Weighill, Yamauchi and Zava~\cite{WeighillYamauchiZava2022}.

\begin{remark}[Banach targets]\label{rem:banach}
The conclusion holds with Hilbert space replaced by any Banach space $X$ of nontrivial type, that is, any $K$-convex space. The metric cotype argument of Mendel and Naor~\cite{MendelNaor2008} evaluates a map only on the finite metric spaces $\{s\sum_j e^{2\pi i x_j/m}e_j\}\subset\ell_\infty^n(\mathbb C)$, which the first assertion of Corollary~\ref{cor:staircase-nonembed} places isometrically in $\mathcal R_1$; on these copies a coarse embedding $f$ of $\mathcal R_1$ into $X$ restricts to maps whose moduli $\omega_f(t)=\inf\{\|f(x)-f(y)\|:d(x,y)\ge t\}$ and $\Omega_f(t)=\sup\{\|f(x)-f(y)\|:d(x,y)\le t\}$ are independent of $n$.
Their Lemma~7.1, with $p=q$ and $r=\infty$, gives $n^{1/q}\omega_f(2s)\le\Gamma m\,\Omega_f(2\pi s/m)$ with $m=m_q(X;n,\Gamma)$ for a constant $\Gamma$, and for every $q$ above the cotype of $X$ one has $n^{1/q}/\Gamma\lesssim m=O(n^{1/q})$ by their Lemma~2.3 and Theorem~4.1, $X$ being $K$-convex.
Taking $s=n^{1/q}$ bounds $\omega_f(2n^{1/q})$ by $C\,\Omega_f(2\pi C'\Gamma)$, with $C,C'$ independent of $n$, so $f$ is not coarse.
\end{remark}

\section{Non-embeddability under a rank condition}\label{sec:rank}

After bars per slice and the size of the presentation, the third complexity bound we consider is a bound on the rank invariant: $\operatorname{rank}M(x\le y)\le n$ for all $x\le y$, which is equivalent to $\dim M(x)\le n$ everywhere (take $x=y$ for one direction, and $\operatorname{rank}M(x\le y)\le\dim M(x)$ for the other).
Unlike a count of bars, it makes sense in any number of parameters.
We show that it is too weak to give a coarse embedding, already in one parameter and already for $n=1$.

In one parameter there are two ways to bound the number of bars by $n$, and they behave differently.
Let $\D{<\infty}$ be the space of persistence diagrams with finitely many off-diagonal points, under the bottleneck distance $\db$; it does not coarsely embed into Hilbert space~\cite{BubenikWagner2020,Mitra2021} (nor, under the $p$-Wasserstein distance with $2<p<\infty$, \cite{Wagner2021}).
Bounding the \emph{total} number of bars gives $\D n$, the diagrams with at most $n$ off-diagonal points, and $\D n$ does coarsely embed~\cite{Mitra2021}.
Bounding the number of bars \emph{alive at any one time} gives a larger space.
For $n=1$, $\D1$ allows a single bar, while the second condition allows any number of bars that never overlap; a diagram with a thousand disjoint bars satisfies it.
Only the second condition makes sense in several parameters, where a module has no bars to count in total; it is the rank bound.
In one parameter it reads $\Rk A(s,t)\le n$ for all $s\le t$, where the \emph{persistent Betti number} $\Rk M(s,t)=\operatorname{rank}(M(s)\to M(t))$ of the module with barcode $A$ is the number of bars $[b,d)$ of $A$ with $b\le s$ and $d>t$.
So we consider
\[
  \Dh{n}:=\{A\in\D{<\infty}:\ \Rk A(s,t)\le n\ \text{for all}\ s\le t\}\ \supseteq\ \D n ,
\]
the finite diagrams with at most $n$ bars alive at any parameter.
The extra room is enough to lose the embedding: diagrams with many disjoint bars, placed far apart, carry every finite metric space, so $\Dh{n}$ contains the family of Dranishnikov, Gong, Lafforgue and Yu and does not coarsely embed into Hilbert space for any $n\ge1$, nor do the two-parameter modules with the same rank bound, even of bounded support (Theorems~\ref{thm:nonembed-diagrams} and~\ref{thm:nonembed-modules}).

The space $\D{<\infty}$ can be described by persistent Betti numbers as well.
Round every coordinate of every diagram down to an integer, deleting the points that land on the diagonal.
This moves each diagram by at most $1$ in $\db$, so the resulting subsets $\D{<\infty}'\subset\D{<\infty}$ and $\D n'\subset\D n$ are coarsely equivalent to the spaces they sit in.
A diagram $A\in\D{<\infty}'$ is determined by the integer table $(\Rk A(s,t))_{s\le t,\ s,t\in\Z}$, because its persistent Betti numbers change only at integer parameters, and the tables that arise have finite entries, almost all zero.
So the non-embeddability of $\D{<\infty}$ says that the one-parameter modules whose persistent Betti numbers are finite and vanish outside a bounded range, the range depending on the module, do not coarsely embed into Hilbert space; neither do the two-parameter modules with the same property, into which the pieces of Theorem~\ref{thm:dgly} embed isometrically and far apart (Theorem~\ref{thm:nonembed-modules}).
If the range is fixed in advance, the space is bounded and embeds trivially; so the lack of a uniform bound on the support matters as much as the unbounded number of bars.

The obstruction we use is the following family of finite metric spaces, which admits no coarse embedding into Hilbert space once its members are placed far apart.

\begin{theorem}[{\cite[Prop.~6.3]{DGLY2002}}] \label{thm:dgly}
Consider $(\mathbb{Z}_n)^m$ as a metric space, where $\mathbb{Z}_n$ carries the word metric and the $m$-fold product the max metric $d_{\infty}$. Let $Y$ be the disjoint union of the $(\mathbb{Z}_n)^m$ (for all $m,n \ge 1$), and let $d$ be a metric on $Y$ whose restriction to each $(\mathbb{Z}_n)^m$ coincides with its existing metric, and such that $d(x,y) \ge m+n+m'+n'$ for $x \in (\mathbb{Z}_n)^m$ and $y \in (\mathbb{Z}_{n'})^{m'}$ with $(n,m)\ne(n',m')$. Then $Y$ does not coarsely embed into Hilbert space.
\end{theorem}

\noindent So a metric space that contains a copy of $Y$, each piece isometrically and the pieces far apart, admits no coarse embedding into Hilbert space either.
Disjoint bars give $\Dh{n}$ room for exactly such a copy.

\begin{theorem}
\label{thm:nonembed-diagrams}
$\Dh{n}$ does not coarsely embed into Hilbert space for any $n \geq 1$.
\end{theorem}

As in \cite[Lemma~4.1]{Mitra2021}, the first step embeds every finite metric space isometrically into $\Dh{1}$, by a variant of the Fr\'echet embedding.

\begin{lemma}\label{lem:finite-isometric}
Every finite metric space $(X,d)$ embeds isometrically into $\Dh{1}$. Furthermore, for any given $a\ge0$ and $b>0$ we can assume that every bar of every diagram in the image of such an embedding has length at least $\diam(X)+b$, and that the persistent Betti numbers $\Rk A(s,t)$ of these diagrams are zero if either $s$ or $t$ lies outside the interval 
\[
\left[a,a+3(\diam(X)+b)|X| + b+2\diam(X)\right].
\]
\end{lemma}

\begin{proof}
Let $X=\{x_1, x_2, \ldots, x_{\kappa}\}$, so that $\kappa=|X|$, and let $R=\diam(X)$. For each $k$ define the persistence diagram
 \[
 f(x_k) = \Big\{
 \big(
 a+3(R+b)\cdot i,
 a+3(R+b)\cdot i + R + b + d(x_k,x_i)
 \big) 
 \mid i=1,2,\ldots,\kappa\Big\}.
 \]
Observe the following two properties of $f$:
\begin{itemize}
    \item The coordinates of $f(x_k)$ (i.e., all the interval endpoints) lie in $[a,a+3(R+b)\kappa + b + 2R]$, and every bar has length $R+b+d(x_k,x_i)\ge R+b$, so the conclusions of the lemma involving $a$ and $b$ hold.
    \item For each $i$, the death coordinate of the point of $f(x_k)$ with birth coordinate $a+3(R+b)\cdot i$ is in the interval
    \[
\mathcal{J}_i=[a+3(R+b)\cdot i + R + b, a+3(R+b)\cdot i + b + 2R]
    \]
    of length $R$, which ends before the next birth coordinate $a+3(R+b)\cdot(i+1)$. So each bar of $f(x_k)$ dies before the next one is born, at most one bar is alive at any parameter, and $f(x_k) \in \Dh{1}$.
\end{itemize}
 
 The map $f \colon X \to f(X)\subset \Dh{1}$ is an isometry, by the following facts:
\begin{enumerate}
 \item For each $k$ the subset $f(x_k)$ consists of precisely one point at each birth coordinate of the form $a+ 3(R+b)\cdot i$, $i=1,2,\ldots,\kappa$, and distinct such birth coordinates are at least $3(R+b)\ge R$ apart.
 \item For each $j,k,i$, 
 \[
\begin{aligned}
&d_\infty\Big(\big(a+3(R+b)\cdot i,\ a+3(R+b)\cdot i + R + b + d(x_j,x_i)\big),\\
&\qquad\big(a+3(R+b)\cdot i,\ a+3(R+b)\cdot i + R + b + d(x_k,x_i)\big)\Big)\\
&\qquad=|d(x_j,x_i)-d(x_k,x_i)| \leq d(x_j,x_k)
\end{aligned}
\] 
 with equality attained at $i=j$ and $i=k$.
 \item A point $(s,t)$ is at bottleneck distance $\tfrac12(t-s)$ from the diagonal, and the point of $f(x_k)$ with birth coordinate $a+3(R+b)\cdot j$ has $t-s=R+b+d(x_k,x_j)$.
\end{enumerate}
By (1) and (2), pairing points of $f(x_j)$ and $f(x_k)$ that share a birth coordinate costs $\max_i|d(x_j,x_i)-d(x_k,x_i)|$, which equals $d(x_j,x_k)$; so $\db(f(x_j),f(x_k))\le d(x_j,x_k)$.
Conversely, in any pairing the point of $f(x_k)$ with birth coordinate $a+3(R+b)\cdot j$ goes either to the point of $f(x_j)$ of the same birth coordinate, at cost $d(x_j,x_k)$ by (2), or to another point of $f(x_j)$, at cost at least $\min\{3(R+b),\tfrac12(R+b+d(x_j,x_k))\}\ge d(x_j,x_k)$ by (1) and (3), since $\db$ caps the $\ell^\infty$-distance by the larger distance to the diagonal, or to the diagonal, at cost $\tfrac12(R+b+d(x_j,x_k))\ge d(x_j,x_k)$ by (3), since $d(x_j,x_k)\le R$.
Hence $d(x_k,x_j)=\db(f(x_j),f(x_k))$ for all $j,k$, and $f$ is an isometry.
\end{proof}

\begin{proof}[Proof of Theorem \ref{thm:nonembed-diagrams}]
It suffices to treat $\Dh{1}$, since $\Dh{1}\subseteq\Dh{n}$ and a coarse embedding restricts to subspaces.

Using the maps of Lemma~\ref{lem:finite-isometric} with suitable $a$ and $b$, we embed isometric copies of $(\mathbb{Z}_n)^m$ for all positive $m,n\in \Z$ into $\Dh{1}$, far from each other and, through $b$, far from the diagonal, so that Theorem~\ref{thm:dgly} applies.
In detail, enumerate the pieces $(\mathbb{Z}_n)^m$ as $P_1,P_2,\ldots$ so that $\nu_i:=n+m$ is nondecreasing in $i$, and let $R_i=\diam(P_i)$ and $b_i=4\nu_i+1$.
Lemma~\ref{lem:finite-isometric} with $b=b_i$ gives an isometry $f_i\colon P_i\to\Dh{1}$ whose bars have length at least $R_i+b_i$ and lie in an interval $I_i$; choosing $a=a_i$ inductively, we may place $I_{i+1}$ to the right of $I_i$ at distance more than $2\nu_{i+1}$, so that $I_i$ and $I_l$ are at distance more than $2\max(\nu_i,\nu_l)$ for $i\ne l$.
Let $x\in P_i$ and $y\in P_l$ with $i\ne l$ and $\nu_i\ge \nu_l$.
In any pairing of $f_i(x)$ and $f_l(y)$, a point of $f_i(x)$ goes either to the diagonal, at cost at least $\tfrac12(R_i+b_i)>2\nu_i$, or to a point of $f_l(y)$, at cost at least the smaller of the distance between $I_i$ and $I_l$ and $\tfrac12(R_i+b_i)$, by the same cap, which exceeds $2\nu_i$.
Hence $\db(f_i(x),f_l(y))>2\nu_i\ge \nu_i+\nu_l$, and $\bigsqcup_if_i(P_i)\subset\Dh{1}$ satisfies the hypotheses of Theorem~\ref{thm:dgly}.
\end{proof}

Theorem~\ref{thm:nonembed-diagrams} transfers to modules, bounded support included.

\begin{theorem}[Modules]\label{thm:nonembed-modules}
Let $\widehat{\mathcal M}_n$ be the finitely presented two-parameter modules $M$ with bounded support and $\Rk M(x,y)\le n$ for all $x\le y$.
Then neither $(\widehat{\mathcal M}_n,\di)$ nor $(\widehat{\mathcal M}_n,\dM)$ coarsely embeds into Hilbert space, for any $n\ge1$.
\end{theorem}

\begin{proof}
For a finite barcode $A$ let $M_A$ be the one-parameter module with barcode $A$, and for $c>0$ let $\widetilde M_{A,c}$ be the two-parameter module with $\widetilde M_{A,c}(x)=M_A(x_1)$ for $0\le x_2<c$ and $0$ otherwise, the structure maps being those of $M_A$ between grades in the strip $\{0\le x_2<c\}$ and zero otherwise.
It is finitely presented, with generators at the $(b,0)$ and relations at the $(d,0)$ and $(b,c)$ for the bars $[b,d)$ of $A$, has bounded support, and $\Rk\widetilde M_{A,c}(x,y)$ equals $\Rk A(x_1,y_1)$ when $0\le x_2\le y_2<c$ and $0$ otherwise; so $\widetilde M_{A,c}\in\widehat{\mathcal M}_n$ when $A\in\Dh{n}$.
Three properties of $A\mapsto\widetilde M_{A,c}$:
\begin{enumerate}
\item $\di(\widetilde M_{A,c},\widetilde M_{B,c})\le\db(A,B)$: an $\eps$-interleaving $(f,g)$ of $M_A$ and $M_B$, which exists for every $\eps>\db(A,B)$ by the isometry theorem~\cite{Lesnick2015}, extends by $\tilde f_x=f_{x_1}$ for $0\le x_2<c-\eps$ and $\tilde f_x=0$ otherwise, and likewise $\tilde g$; naturality and the interleaving identities hold grade by grade, the vertical structure maps of the strip being identities or zero and the horizontal ones those of $M_A$.
\item If $\db(A,B)<c/2$ then equality holds: by (1), $\eps$-interleavings exist for every $\eps\in(\di,c/2)$, and an $\eps$-interleaving $(\tilde f,\tilde g)$ of $\widetilde M_{A,c}$ and $\widetilde M_{B,c}$ with $\eps<c/2$ restricts along the line $x_2=0$ to an $\eps$-interleaving of $M_A$ and $M_B$, since $\widetilde M_{B,c}[\eps](t,0)=M_B(t+\eps)$ and $\widetilde M_{A,c}[2\eps](t,0)=M_A(t+2\eps)$, and naturality along $(t,0)\le(t,\eps)$ identifies $\tilde g_{(t,\eps)}$ with $\tilde g_{(t,0)}$; so $\db(A,B)\le\di(\widetilde M_{A,c},\widetilde M_{B,c})$ by the isometry theorem.
\item If the bars of $A$ and of $B$ lie in intervals at distance $\delta>\eps$, an $\eps$-interleaving of $\widetilde M_{A,c}$ and $\widetilde M_{B,c'}$ has $\tilde f=\tilde g=0$, the supports being disjoint after the shift, so the $2\eps$-shift maps of both modules vanish; but the shift map of $\widetilde M_{A,c}$ is nonzero as long as $2\eps<c$ and $2\eps$ is less than the length $\lambda_A$ of the longest bar of $A$.
Hence $\di(\widetilde M_{A,c},\widetilde M_{B,c'})\ge\min\bigl\{\delta,\ \tfrac12\max(\min(c,\lambda_A),\min(c',\lambda_B))\bigr\}$.
\end{enumerate}
Now take the pieces $f_i(P_i)\subset\Dh{1}$ of the proof of Theorem~\ref{thm:nonembed-diagrams}, whose bars have length at least $R_i+b_i$, put $c_i:=2R_i+4\nu_i+1$, and let $F_i(x):=\widetilde M_{f_i(x),c_i}$.
Within a piece, $\db(f_i(x),f_i(y))=d(x,y)\le R_i<c_i/2$, so $F_i$ is an isometry by (1) and (2).
Across pieces $P_i,P_l$ with $\nu_i\ge \nu_l$, the intervals containing the bars are at distance more than $2\nu_i$ and $\tfrac12\min(c_i,R_i+b_i)=\tfrac12\min(2R_i+4\nu_i+1,R_i+4\nu_i+1)>2\nu_i$, so $\di(F_i(x),F_l(y))>2\nu_i\ge \nu_i+\nu_l$ by (3).
Thus $\bigsqcup_iF_i(P_i)\subset\widehat{\mathcal M}_n$ with the metric $\di$ satisfies the hypotheses of Theorem~\ref{thm:dgly}, and $(\widehat{\mathcal M}_n,\di)$ does not coarsely embed into Hilbert space.

For $\dM$, which is at most $\di$ (Definition~\ref{def:dmatch}), it suffices to bound it from below on the same family.
Let $\beta$ be the birth of the $j$-th bar of the diagrams in $f_i(P_i)$ and $L(t)=(\beta,0)+t\one$, a normalized line.
Along $L$ the module $F_i(x)$ restricts to the $j$-th bar of $f_i(x)$ alone, as $[0,\lambda)$ with $\lambda$ its length: that bar has length at most $2R_i+b_i=c_i$, so it dies inside the strip, the earlier bars die before $\beta$, and the next is born $3(R_i+b_i)>c_i$ later.
For $x=x_k$, $y=x_{k'}$ in $P_i$ and $j=k$ the two restrictions are $[0,R_i+b_i)$ and $[0,R_i+b_i+d(x,y))$, at bottleneck distance $\min\{d(x,y),\tfrac12(R_i+b_i+d(x,y))\}=d(x,y)$, so $\dM(F_i(x),F_i(y))=d(x,y)$.
For $y\in P_l$, $l\ne i$, $\nu_i\ge \nu_l$, the bars of $F_l(y)$ along $L$ start at parameters $t\ge0$ with $\beta+t\in I_l$: there are none if $I_l$ lies to the left of $I_i$, and otherwise they start at $t>2\nu_i$; the bar $[0,\lambda)$ of $F_i(x)$ has $\lambda\ge R_i+b_i>4\nu_i$, so matching it to the diagonal or to a bar of $F_l(y)$ costs more than $2\nu_i$; hence $\dM(F_i(x),F_l(y))>2\nu_i$, and Theorem~\ref{thm:dgly} applies again.
\end{proof}

\section{Discussion}\label{sec:disc}
How complexity is bounded decides the coarse geometry of multiparameter persistence.
On every finite net of slices, a bound on the presentation, or merely on bars per slice, gives a quantitative coarse embedding into Hilbert space and finite asymptotic dimension (Theorem~\ref{thm:modules}).
Under the continuous matching and interleaving distances, a bound on bars per slice gives neither: the staircases have one bar per slice, infinite asymptotic dimension and no coarse embedding (Theorem~\ref{thm:staircase}, Corollary~\ref{cor:staircase-nonembed}).
A bound on persistent Betti numbers gives no coarse embedding either, even for modules of bounded support (Theorems~\ref{thm:nonembed-diagrams} and~\ref{thm:nonembed-modules}).

The embeddings certify separation only in the sliced distance $\dM^S\le\dM$, and their lower gauge carries the factor $3^{-N/2}$, which decays with the complexity $N$.
Whether a bounded presentation gives finite asymptotic dimension under the continuous matching distance is open; the companion paper~\cite{MajhiMitraVirk-asdim} treats it together with the interleaving distance.

\bibliographystyle{plain}
\bibliography{main}

\begin{thebibliography}{10}

\bibitem{Assouad1983}
Patrice Assouad.
\newblock Plongements lipschitziens dans $\mathbb{R}^n$.
\newblock {\em Bulletin de la Soci{\'e}t{\'e} Math{\'e}matique de France},
  111:429--448, 1983.

\bibitem{BateGarciaPulido2024}
David Bate and Ana~Luc{\'i}a Garcia~Pulido.
\newblock Bi-{L}ipschitz embeddings of the space of unordered $m$-tuples with a
  partial transportation metric.
\newblock {\em Mathematische Annalen}, 390(2):3109--3131, 2024.

\bibitem{BellDranishnikov2008}
G.~Bell and A.~Dranishnikov.
\newblock Asymptotic dimension.
\newblock {\em Topology and its Applications}, 155(12):1265--1296, 2008.

\bibitem{BLPY2019}
Greg Bell, Austin Lawson, C.~Neil Pritchard, and Dan Yasaki.
\newblock The space of persistence diagrams fails to have {Y}u's property {A}.
\newblock {\em Topology Proceedings}, 58:279--288, 2021.
\newblock arXiv:1902.02288.

\bibitem{BjerkevikBotnanKerber2020}
H\aa{}vard~Bakke Bjerkevik, Magnus~Bakke Botnan, and Michael Kerber.
\newblock Computing the interleaving distance is {NP}-hard.
\newblock {\em Foundations of Computational Mathematics}, 20(5):1237--1271,
  2020.

\bibitem{BjerkevikLesnick2021}
H{\aa}vard~Bakke Bjerkevik and Michael Lesnick.
\newblock $\ell^p$-distances on multiparameter persistence modules.
\newblock arXiv:2106.13589, 2021.

\bibitem{BotnanLesnick2023}
Magnus~Bakke Botnan and Michael Lesnick.
\newblock An introduction to multiparameter persistence.
\newblock In Aslak~Bakke Buan, Henning Krause, and {\O}yvind Solberg, editors,
  {\em Representations of Algebras and Related Structures}, EMS Series of
  Congress Reports, pages 77--150. EMS Press, Berlin, 2023.
\newblock arXiv:2203.14289.

\bibitem{Botnan2022signed}
Magnus~Bakke Botnan, Steffen Oppermann, and Steve Oudot.
\newblock Signed barcodes for multi-parameter persistence via rank
  decompositions.
\newblock In {\em 38th International Symposium on Computational Geometry (SoCG
  2022)}, volume 224 of {\em Leibniz International Proceedings in Informatics
  (LIPIcs)}, pages 19:1--19:18. Schloss Dagstuhl--Leibniz-Zentrum f\"ur
  Informatik, 2022.

\bibitem{BDLM2008}
N.~Brodskiy, J.~Dydak, M.~Levin, and A.~Mitra.
\newblock A {H}urewicz theorem for the {A}ssouad--{N}agata dimension.
\newblock {\em Journal of the London Mathematical Society}, 77(3):741--756,
  2008.

\bibitem{BubenikDeSilvaScott2015}
Peter Bubenik, Vin de~Silva, and Jonathan Scott.
\newblock Metrics for generalized persistence modules.
\newblock {\em Foundations of Computational Mathematics}, 15(6):1501--1531,
  2015.

\bibitem{BubenikWagner2020}
Peter Bubenik and Alexander Wagner.
\newblock Embeddings of persistence diagrams into {H}ilbert spaces.
\newblock {\em Journal of Applied and Computational Topology}, 4(3):339--351,
  2020.

\bibitem{CarlssonZomorodian2009-multidim}
Gunnar Carlsson and Afra Zomorodian.
\newblock The theory of multidimensional persistence.
\newblock {\em Discrete \& Computational Geometry}, 42(1):71--93, 2009.

\bibitem{Carriere2019-nx}
Mathieu Carri{\`e}re and Ulrich Bauer.
\newblock On the metric distortion of embedding persistence diagrams into
  separable {H}ilbert spaces.
\newblock In Gill Barequet and Yusu Wang, editors, {\em 35th International
  Symposium on Computational Geometry (SoCG 2019)}, volume 129 of {\em LIPIcs},
  pages 21:1--21:15. Schloss Dagstuhl--Leibniz-Zentrum f{\"u}r Informatik,
  2019.

\bibitem{Carriere2020-mpi}
Mathieu Carri\`ere and Andrew~J. Blumberg.
\newblock Multiparameter persistence image for topological machine learning.
\newblock In {\em Advances in Neural Information Processing Systems (NeurIPS)},
  volume~33, 2020.

\bibitem{CerriDiFabioFerriEtAl2013}
Andrea Cerri, Barbara Di~Fabio, Massimo Ferri, Patrizio Frosini, and Claudia
  Landi.
\newblock Betti numbers in multidimensional persistent homology are stable
  functions.
\newblock {\em Mathematical Methods in the Applied Sciences},
  36(12):1543--1557, 2013.

\bibitem{DGLY2002}
A.~N. Dranishnikov, G.~Gong, V.~Lafforgue, and G.~Yu.
\newblock Uniform embeddings into {H}ilbert space and a question of {G}romov.
\newblock {\em Canadian Mathematical Bulletin}, 45(1):60--70, 2002.

\bibitem{Gal2008}
{\'S}wiatos{\l}aw~R. Gal.
\newblock Asymptotic dimension and uniform embeddings.
\newblock {\em Groups, Geometry, and Dynamics}, 2(1):63--84, 2008.

\bibitem{Gonzalez1985}
Teofilo~F. Gonzalez.
\newblock Clustering to minimize the maximum intercluster distance.
\newblock {\em Theoretical Computer Science}, 38:293--306, 1985.

\bibitem{Heinonen2001}
Juha Heinonen.
\newblock {\em Lectures on Analysis on Metric Spaces}.
\newblock Universitext. Springer-Verlag, New York, 2001.

\bibitem{HurewiczWallman1941}
Witold Hurewicz and Henry Wallman.
\newblock {\em Dimension Theory}, volume~4 of {\em Princeton Mathematical
  Series}.
\newblock Princeton University Press, Princeton, NJ, 1941.

\bibitem{KerberLesnickOudot2019}
Michael Kerber, Michael Lesnick, and Steve Oudot.
\newblock Exact computation of the matching distance on 2-parameter persistence
  modules.
\newblock In {\em 35th International Symposium on Computational Geometry (SoCG
  2019)}, volume 129 of {\em Leibniz International Proceedings in Informatics
  (LIPIcs)}, pages 46:1--46:15. Schloss Dagstuhl--Leibniz-Zentrum f\"ur
  Informatik, 2019.

\bibitem{Korkmaz2025-cumperlay}
Caner Korkmaz, Brighton Nuwagira, Bar{\i}{\c{s}} Co{\c{s}}kunuzer, and Tolga
  Birdal.
\newblock {CuMPerLay}: Learning cubical multiparameter persistence
  vectorizations.
\newblock arXiv:2510.12795, 2025.

\bibitem{Landi2018}
Claudia Landi.
\newblock The rank invariant stability via interleavings.
\newblock In {\em Research in Computational Topology}, volume~13 of {\em
  Association for Women in Mathematics Series}, pages 1--10. Springer, 2018.

\bibitem{LangSchlichenmaier2005}
Urs Lang and Thilo Schlichenmaier.
\newblock Nagata dimension, quasisymmetric embeddings, and {L}ipschitz
  extensions.
\newblock {\em International Mathematics Research Notices},
  2005(58):3625--3655, 2005.

\bibitem{Lesnick2015}
Michael Lesnick.
\newblock The theory of the interleaving distance on multidimensional
  persistence modules.
\newblock {\em Foundations of Computational Mathematics}, 15(3):613--650, 2015.

\bibitem{Loiseaux2023}
David Loiseaux, Luis Scoccola, Mathieu Carri\`ere, Magnus~Bakke Botnan, and
  Steve Oudot.
\newblock Stable vectorization of multiparameter persistent homology using
  signed barcodes as measures.
\newblock In {\em Advances in Neural Information Processing Systems (NeurIPS)},
  volume~36, 2023.

\bibitem{MajhiMitraVirk-asdim}
Sushovan Majhi, Atish Mitra, and {\v{Z}}iga Virk.
\newblock Asymptotic dimension of persistence modules of bounded presentation.
\newblock arXiv:2611.16965, 2026.

\bibitem{PaperII}
Sushovan Majhi, Atish Mitra, {\v Z}iga Virk, and Pramita Bagchi.
\newblock A closed-form adaptive-landmark kernel for certified point-cloud and
  graph classification.
\newblock arXiv:2605.04046, 2026.

\bibitem{PaperI}
Sushovan Majhi, Atish Mitra, {\v Z}iga Virk, and Pramita Bagchi.
\newblock A closed-form persistence-landmark pipeline for certified point-cloud
  and graph classification.
\newblock {\em Transactions on Machine Learning Research}, 2026.

\bibitem{Majhi2026-complex}
Sushovan Majhi, Atish Mitra, {\v Z}iga Virk, and Pramita Bagchi.
\newblock {COMPLEX}: A closed-form certified embedding of multiparameter
  persistence modules.
\newblock arXiv:2609.22012, 2026.

\bibitem{MendelNaor2008}
Manor Mendel and Assaf Naor.
\newblock Metric cotype.
\newblock {\em Annals of Mathematics}, 168(1):247--298, 2008.

\bibitem{MillerSturmfels2005}
Ezra Miller and Bernd Sturmfels.
\newblock {\em Combinatorial Commutative Algebra}, volume 227 of {\em Graduate
  Texts in Mathematics}.
\newblock Springer, New York, 2005.

\bibitem{Mitra2021}
Atish Mitra and {\v Z}iga Virk.
\newblock The space of persistence diagrams on $n$ points coarsely embeds into
  {H}ilbert space.
\newblock {\em Proceedings of the American Mathematical Society},
  149(6):2693--2703, 2021.

\bibitem{MitraVirk2024corrigendum}
Atish Mitra and {\v Z}iga Virk.
\newblock Corrigendum to ``{T}he space of persistence diagrams on $n$ points
  coarsely embeds into {H}ilbert space''.
\newblock {\em Proceedings of the American Mathematical Society},
  152(4):1803--1807, 2024.

\bibitem{Mitra2024}
Atish Mitra and {\v Z}iga Virk.
\newblock Geometric embeddings of spaces of persistence diagrams with explicit
  distortions.
\newblock arXiv:2401.05298v2, 2024.

\bibitem{Mukherjee2024}
Soham Mukherjee, Shreyas~N. Samaga, Cheng Xin, Steve Oudot, and Tamal~K. Dey.
\newblock {D-GRIL}: End-to-end topological learning with 2-parameter
  persistence.
\newblock In {\em 42nd International Symposium on Computational Geometry
  (SoCG)}, volume 367 of {\em Leibniz International Proceedings in Informatics
  (LIPIcs)}, pages 79:1--79:17. Schloss Dagstuhl -- Leibniz-Zentrum f{\"u}r
  Informatik, 2026.
\newblock arXiv:2406.07100.

\bibitem{NaorNeiman2012}
Assaf Naor and Ofer Neiman.
\newblock Assouad's theorem with dimension independent of the snowflaking.
\newblock {\em Revista Matem{\'a}tica Iberoamericana}, 28(4):1123--1142, 2012.

\bibitem{Roe2003}
John Roe.
\newblock {\em Lectures on Coarse Geometry}, volume~31 of {\em University
  Lecture Series}.
\newblock American Mathematical Society, Providence, RI, 2003.

\bibitem{Vipond2020}
Oliver Vipond.
\newblock Multiparameter persistence landscapes.
\newblock {\em Journal of Machine Learning Research}, 21(61):1--38, 2020.

\bibitem{Vipond2021}
Oliver Vipond, Joshua~A. Bull, Philip~S. Macklin, Ulrike Tillmann,
  Christopher~W. Pugh, Helen~M. Byrne, and Heather~A. Harrington.
\newblock Multiparameter persistent homology landscapes identify immune cell
  spatial patterns in tumors.
\newblock {\em Proceedings of the National Academy of Sciences},
  118(41):e2102166118, 2021.

\bibitem{Wagner2021}
Alexander Wagner.
\newblock Nonembeddability of persistence diagrams with {$p>2$} {W}asserstein
  metric.
\newblock {\em Proceedings of the American Mathematical Society},
  149(6):2673--2677, 2021.

\bibitem{WeighillYamauchiZava2022}
Thomas Weighill, Takamitsu Yamauchi, and Nicol{\`o} Zava.
\newblock Coarse infinite-dimensionality of hyperspaces of finite subsets.
\newblock {\em European Journal of Mathematics}, 8(1):335--355, 2022.

\bibitem{Xin2023GRIL}
Cheng Xin, Soham Mukherjee, Shreyas~N. Samaga, and Tamal~K. Dey.
\newblock {GRIL}: A 2-parameter persistence based vectorization for machine
  learning.
\newblock In {\em Proceedings of the 2nd Annual Workshop on Topology, Algebra,
  and Geometry in Machine Learning (TAG-ML)}, volume 221 of {\em Proceedings of
  Machine Learning Research (PMLR)}, pages 313--333, 2023.

\bibitem{Zava2025}
Nicol{\`o} Zava.
\newblock Coarse and bi-{L}ipschitz embeddability of subspaces of the
  {G}romov--{H}ausdorff space into {H}ilbert spaces.
\newblock {\em Algebraic \& Geometric Topology}, 25(8):5153--5174, 2025.

\end{thebibliography}

\end{document}